\documentclass[12pt]{amsart} 
\usepackage{amsmath, amstext, amsbsy, amssymb}

\newtheorem{theorem}{Theorem}[section]
\newtheorem{lemma}[theorem]{Lemma}

\newtheorem{proposition}{Proposition}[section]
\newtheorem{corollary}[theorem]{Corollary}

\theoremstyle{definition}

\theoremstyle{remark}
\newtheorem{remark}[theorem]{Remark}

\numberwithin{equation}{section}

\newcommand{\C}{\mathbb C}
\newcommand{\ch}{\mbox{ch}}

\newcommand{\bG}{H\G}
\newcommand{\bFG}{{\mathcal F}_{\bG}^T}
\newcommand{\FGG}{\overline{\mathcal F}_{\bG} }
\newcommand{\bFGG}{\overline{\mathcal F}_{\bG}^T }

\newcommand{\G}{\Gamma}
\newcommand{\Gbar}{\overline{\Gamma}^*}

\newcommand{\Gn}{{\Gamma}_n}

\newcommand{\g}{\gamma}
\newcommand{\ga}{\gamma}
\newcommand{\thg}{\widehat{\mathfrak h}_{\G, \wt}[-1]}
\newcommand{\hg}{\widehat{\mathfrak h}_{\G, \wt}}

\newcommand{\thhg}{\widehat{\widehat{\mathfrak g}}[-1]}
\newcommand{\la}{\lambda}

\newcommand{\tloopg}{\widehat{\mathfrak g}[-1]}

\newcommand{\be}{\beta}

\newcommand{\ep}{\epsilon}
\newcommand{\RG}{R_{\bG}}
\newcommand{\tRG}{R^-_{H\Gamma}}  
\newcommand{\RGG}{\overline{R}_{ \bG}}
\newcommand{\tRGG}{\overline{R}_{\bG}^T}
\newcommand{\RGO}{R^0_{ \bG}}
\newcommand{\Rz}{R_{\mathbb Z} (\G)}
\newcommand{\Rzz}{\overline{R}_{\mathbb Z} (\G)}

\newcommand{\tSG}{ S^-_{\bG} }

\newcommand{\tSGG}{\overline{S}_{\bG} }
\newcommand{\SGO}{ S^0_{\bG} }

\newcommand{\tVG}{V_{\bG}^T}

\newcommand{\tVGG}{\overline{V}_{ \bG}^T}

\newcommand{\wt}{\xi}
\newcommand{\Z}{\mathbb Z}

\def\tS#1{\widetilde{S}_{#1}}      
\def\tG#1{\widetilde{\Gamma}_{#1}}
\newcommand{\htimes}{\hat{\times}}
\newcommand{\Rtz}{R_{\mathbb F_2} (\G)}

\newcommand{\Rtzh}{\hat{R}_{\mathbb F_2}^-(\G)}
\newcommand{\Rtzhbar}{\hat{R}_{\mathbb F_2}^-(\Gbar)}
\newcommand{\oa}{\overline{\alpha}}
\newcommand{\ob}{\overline{\beta}}
\def\B#1{\widetilde{H}_{#1}}
\def\BG#1{\widetilde{H}\G_{#1}}
\def\hgr#1{H_{#1}}
\def\HG#1{H\G_{#1}}

\begin{document}

\title[Twisted vertex representations and spin characters]
    {Twisted vertex representations and spin characters}
\author{Naihuan Jing}
\address{Department of Mathematics,
   North Carolina State Univer\-sity,
   Ra\-leigh, NC 27695-8205}
\email{jing@math.ncsu.edu, wqwang@math.ncsu.edu}
\thanks{Jing acknowledges the support by NSF grant
DMS-9970493.}
\author{Weiqiang Wang}

\begin{abstract}
We establish a new group-theoretic realization of
the basic representations of the twisted affine and twisted toroidal algebras
of ADE types in the same spirit of our new approach 
to the McKay correspondence. 
Our vertex operator construction provides
a unified description to the character tables for
the spin cover of the wreath product of
the twisted hyperoctahedral groups and an arbitrary finite group.
\end{abstract}

\maketitle

\section*{Introduction} \label{S:intro}
The connection between the homogeneous vertex representations of
affine Lie algebras and the wreath products $\Gn = \G^n\rtimes
S_n$ associated to finite subgroups $\G$ of $SL_2(\C)$ was first
pointed out in \cite{W} and subsequently established fully in
\cite{FJW1}. This initiates a new approach to the McKay
correspondence (cf. \cite{Mc}), which classically relates the
finite subgroups of $SL_2(\C)$ in a bijective manner with affine
Dynkin diagrams of ADE type.

These results have been further extended in \cite{FJW2} to realize
the vertex representations of twisted affine Lie algebras
$\tloopg$ and its toroidal counterpart by using the spin
representations of a double cover $\tG n = \G^n\rtimes\tS n  $ of
the wreath product $\Gn$. Here $\tS n$ is a double cover of the
symmetric group $S_n$ whose representation theory was developed by
Schur \cite{Sc} and reformulated in terms of twisted vertex
operators in \cite{J}. The {\em algebraic} construction of the
vertex representation of the twisted affine algebra $\tloopg$ was
obtained first in \cite{LW} in the case of $\mathfrak g = sl_2$, and
in \cite{FLM1, FLM2} for general $\mathfrak g$ as one ingredient in
the vertex (operator) algebra construction of the Monster group.

The goal of this paper is to provide a new finite-group-theoretic
realization of the vertex representations of the twisted affine
Lie algebra $\tloopg$ and its toroidal counterpart. Instead of
$\tG n$ we will use a semi-direct product $\BG n =\G^n\rtimes \B
n $, where $\B n$ is a double cover of the hyperoctahedral group
$H_n$. The finite group $\BG n$ can also be thought as a double
cover of the wreath product $(\G \times \Z_2 )^n\rtimes S_n$.

Our present construction recaptures all the constructions in
\cite{FJW2} with additional advantages. First, this has a natural
generalization given in a companion paper \cite{W2}, which is
analogous to \cite{W}, to the $K$-theory setup and it is
intimately related to geometry. Secondly, the present
constructions of the twisted vertex operators, Heisenberg algebra,
and the characteristic map etc are essentially over $\Z$ while
those in \cite{FJW2} involve an inevitable square root of $2$
which originates in the theory of spin representations of $\tS n$.
Our present construction in principle paves the way to study the
quantum twisted vertex representations at roots of unity which we
will investigate elsewhere.

More explicitly, we take the supermodule approach of Joz\'efiak
\cite{Jo1} to study the spin representations of $\BG n$ for an
arbitrary finite group $\G$. We give a description with a complete
proof of the {\em split} conjugacy classes in $\BG n$ which play
an important role in understanding the spin supermodules of $\BG
n$. This generalizes earlier works on $\B n$ and its spin
representations (cf. \cite{Sg, Jo2, St}). This result was also
stated in Read \cite{R} when $\G$ is cyclic.

Given a finite group $\G$, we consider a direct sum over all $n$,
denoted by $\tRG$, of the Grothendieck groups of {\em spin}
supermodules of $\BG n$. We show that $\tRG$ carries a natural
Hopf algebra structure. Associated to a self-dual virtual
character $\xi$ of $\G$ we define a $\xi$-weighted symmetric
bilinear form on $\tRG$ (compare \cite{FJW1, FJW2}). The space
$\tRG$ can be shown to be isomorphic to a Fock space of a twisted
Heisenberg algebra. The twisted vertex operators also make a
natural appearance. One can further interpret the Fock space as a
distinguished space of symmetric functions parametrized by the set
of irreducible characters of $\G$, and in this way the irreducible
characters of $\BG n$ correspond essentially to the Schur
$Q$-functions.

When $\G$ is a subgroup of $SL_2(\C)$ and the weight $\xi$ is
chosen suitably, the Fock space $\tRG$ leads to a construction of
the twisted vertex representations of $\tloopg$ and its toroidal
counterpart. On the other hand, when choosing the weight $\xi$ to
be trivial, the twisted vertex operator approach allows us to
compute spin characters of all irreducible supermodules of the
group $\BG n$, as done for $\tS n$ in \cite{J} and for $\tG n$ in
\cite{FJW2}.

The layout of this paper is as follows. In Sect.~\ref{S:hyperoct}
we studied in detail the conjugacy classes of the group $\BG n$.
In Sect.~\ref{sec_spinsuper}, we introduce the Grothendieck group
$R^-(\BG n)$ and a weighted bilinear form on it. We construct the
Hopf algebra structure on $\tRG$. In Sect.~\ref{sec_heis} we
identify $\tRG$ with a Fock space of a twisted Heisenberg algebra.
In Sect.~\ref{sec_vertex} we construct twisted vertex operators
(essentially on $\tRG$) in terms of group theoretic operators.
When $\G$ is a finite subgroup of $SL_2(\C)$, this leads to a
group theoretic construction of the vertex representation of
twisted affine and toroidal Lie algebras. In Sect.~\ref{sec_char}
we construct the irreducible spin super characters of $\BG n$ and
recover their character table from vertex operator viewpoint. When
the statements can be proved similarly as in \cite{FJW2}, we often
sketch only or omit the proofs and refer the reader to {\em loc.
cit.} for more detail.
\section{The conjugacy classes of the finite groups $\BG n$}
\label{S:hyperoct}

In this section we introduce the finite group $\BG n$ and a
natrual $\Z_2$-grading on it. We also classify the so-called {\em
split} conjugacy classes in $\BG n$ which will play a key role in
the study of spin supermodules of $\BG n$ in later sections.
\subsection{Definitions}
Let $\Pi_n$ be the finite group generated by $a_i$ ($i=1, \ldots,
n$) and the central element $z$ subject to the relations
\begin{equation}\label{E:rel1}
a_i^2=z, \quad z^2=1, \quad a_ia_j=za_ja_i \quad (i\neq j).
\end{equation}
The symmetric group $S_n$ acts on $\Pi_n$ by $s(a_i)=a_{s(i)}$,
$s\in S_n$. The semidirect product $\B n:=\Pi_n\rtimes S_n $ is
called the {\it twisted hyperoctahedral group}. Explicitly the
multiplication in $\B n$ is given by
\begin{equation*}
(a, s)(a', s')=(as(a'), ss'), \qquad a, a'\in \Pi_n, s, s'\in S_n.
\end{equation*}
Since $\Pi_n/\{1, z \} \simeq \Z_2^n$, the group $\B n$ is a
double cover of the hyperoctahedral group $\hgr n:=\Z_2^n\rtimes S_n$.

Let $\G$ be a finite group. The twisted hyperoctahedral group $\B
n$ acts on the product group $\G^n:=\G\times \cdots\times \G$ by
letting $\Pi_n$ act trivially on $\G^n$ and letting $S_n$ act by
$s(g_1, \cdots, g_n) =(g_{s^{-1}(1)}, \cdots, g_{s^{-1}(n)}), s\in
S_n.$ The finite group $\BG n$ is then defined to be the
semi-direct product of  $\B n$ and $\G^n$. Alternatively, the
symmetric group $S_n$ naturally acts on $\G^n \times \Pi_n $ by
simutaneous permutations of elements in $\G^n$ and $ \Pi_n $, and
we may regard $\BG n$ as the semi-direct product of the symmetric
group $S_n$ and $\G^n \times \Pi_n $.

The double covering $\B n$ of $\hgr n$ extends to a double
covering of the wreath product $\HG n :=(\G\times
\mathbb Z_2)^n\rtimes S_n$ by $\BG n$:
\begin{equation*}
1 \longrightarrow\mathbb Z_2 \stackrel{i}{\longrightarrow} \BG n
\stackrel{\theta_n}{\longrightarrow} \HG n \longrightarrow 1.
\end{equation*}
The order $|\B n|$ is clearly $2^{n+1}n! |\G|^n$, where $|\G|$
denotes the order of $\G$. The group $\BG n$ contains several
distinguished subgroups: $S_n$, $\B n$, $\G^n$, $\BG n$, and the
wreath product $\Gn :=\G^n\rtimes S_n$ etc.

We define a $\mathbb Z_2$-grading on the group $\BG n$ by setting
the degree of $a_i$ to be $1$ and the degree of elements in $S_n$
and $\G^n$ to be $ 0$, and denote by $\BG n^0$ (resp. $\BG n^1$)
the degree zero (resp. one) part. This induces an epimorphism $p$
>from $\BG n$ or the group algebra $\mathbb C[\BG n]$ to $\mathbb
Z_2$. The homomorphism $p$ descends to a homomorphism on $\HG n$
which will be denoted by $p$ again. We say $x\in \BG n$ or $\HG n$
{\it even} (resp. {\it odd}) if $p(x)$ is $0$ (resp. $1$).
\subsection{Partition-valued functions}
A partition $\la$ of a non-negative integer $n$ is a monotonic
non-increasing sequence of integers $\la_i$ called parts such that
$n=\la_1+\cdots+\la_l=|\la|$. Here $l=l(\la)$ is the {\it length}
of $\la$. We may also write $\la=(1^{m_1}2^{m_2}\cdots)$, where
$m_i$ is the number of times that $i$ appears in $\la$. For two
partitions $\la$ and $\mu$ the dominance order $\la\geqslant\mu$
is defined by $\la_1\geq \mu_1$, $\la_1+\la_2\geq \mu_1+\mu_2$,
etc. A partition $\la$ is called {\it strict} if its parts are
distinct integers.

We will use partitions indexed by $\G_*$ and $\G^*$. For a finite
set $X$ and $\rho=(\rho(x))_{x\in X}$ a family of partitions
indexed by $X$, we write $\|\rho\|=\sum_{x\in X}|\rho(x)|.$ It is
convenient to regard $\rho=(\rho(x))_{x\in X}$ as a
partition-valued function on $X$. We denote by $\mathcal{P}(X)$
the set of all partitions indexed by $X$ and by $\mathcal{P}_n(X)$
the subset consisting of $\rho$ such that $\|\rho\|=n$. The total
number of parts, denoted by $l(\rho)=\sum_xl(\rho(x))$, in the
partition-valued function $\rho=(\rho(x))_{x\in X}$ is called the
{\it length} of $\rho$. The {\it dominance order} on
$\mathcal{P}(X)$ is defined naturally by $\rho\geqslant \pi$ if
$\rho(x)\geqslant \pi(x)$ for each $x$. We also write $\rho\gg\pi$
if $\rho(x)\geqslant \pi(x)$ and $\rho(x)\neq\pi(x)$ for each
$x\in X$. For $\rho\in \mathcal {OP} ( \G_* )$ we define
$\overline{\rho}  \in  \mathcal {OP} ( \G_* )$ by
$\overline{\rho}(c)=\rho (c^{-1})$, where $c^{-1}$ denotes the
conjugacy class $\{g|g^{-1} \in c \}$.

Let $\mathcal{OP}(X)$ be the set of partition-valued functions
$(\rho(x))_{x\in X}$ in $\mathcal{P}(X)$ such that all parts of
$\rho(x)$ are odd integers for each $x$, and let $\mathcal{SP}(X)$
be the set of $\rho$ such that each $\rho(x)$ is strict. When $X$
consists of a single element, we will omit $X$ and simply write
$\mathcal P$ for $\mathcal P(X)$, thus $\mathcal OP$ or $\mathcal
SP$ will be used accordingly. A variant of Euler's theorem says
that the number of strict partition-valued functions on a set $X$
is equal to the number of partition-valued functions on $X$ with
odd integer parts.

We denote by
\begin{eqnarray*}
{\mathcal P}_n^+(X)&=\{\la\in {\mathcal P}_n(X)|\quad
l(\rho)\equiv 0 \bmod{2} \}, \\
 {\mathcal P}_n^-(X) &=\{\la\in {\mathcal P}_n(X)|\quad l(\rho)\equiv 1
\bmod{2} \},
\end{eqnarray*}
and define ${\mathcal SP}_n^{\pm}(X)={\mathcal P}_n^{\pm}(X)\cap
{\mathcal SP}_n(X)$ for $i=0, 1$.

We will also need another parity $d$ on partition-valued
functions. For a partition-valued function $\rho=(\rho(x))_{x\in
X}$ we define
\begin{equation}\label{E:parity1}
  d(\rho)=\|\rho\|-l(\rho).
\end{equation}
\subsection{Conjugacy classes of $\HG n$} \label{S:conj}

Let $\Gamma$ be a finite group with $r +1$ conjugacy classes. We
denote by $\Gamma^*=\{\g_i\}_{i=0}^{r}$ the set of complex
irreducible characters where $\g_0$ is the trivial character, and
by $\Gamma_*=\{c^i\}$ the set of conjugacy classes where $c^0$ is
the identity conjugacy class. Let $|c|$ be the order of the
conjugacy class $c\in\G_*$, and then $\zeta_c :=|\G |/|c|$ is the
order of the centralizer of an element in the class $c$.

For a subset $I= \{i_1,i_2, \ldots, i_m\}$ of the set $\{1,
\cdots, n\}$, we denote $a_I =a_{i_1}a_{i_2}\cdots a_{i_m}$. It
follows that $p(a_I)\equiv |I|\, (mod\, 2)$. If $I\cap J=\emptyset
$, then $a_Ia_J=(-1)^{|I||J|}a_Ja_I$. Also we can easily show by
induction that

\begin{eqnarray}  \label{eq_perm}
a_{i_1 i_2\cdots i_m}=z^{d(s)} a_{s(i_1)s(i_2)\cdots s(i_m)}
\end{eqnarray}
for a permutation $s$ of $I=
\{i_1,i_2, \ldots, i_m\}$. Clearly the morphism $\theta_n :\B n
\rightarrow H_n$ sends $a_I$ to $b_I=b_{i_1}\cdots b_{i_m}$, where
$b_i$ are the generators of $\mathbb Z_2^n$.

The conjugacy classes of a wreath product is well understood, cf.
\cite{M, Z}. In particular this gives us the following description
of conjugacy classes of the wreath product $\HG n =(\G
\times \Z_2)^n\rtimes S_n$. Given a cycle $t= (i_1 \ldots i_m)$, we call the
set $\{i_1, \ldots, i_m \}$ the {\em support} of $t$, denoted by
$supp(t)$. Given $(g, b_Is)\in \HG n$, every element $b_Is \in H_n$
can be uniquely written as a product (up to order)
\[ b_Is= (b_{I_1}s_1)\cdots (b_{I_k}s_k),
\]
where $s \in S_n$ is a product of disjoint cycles $s_1 \ldots
s_k$, and $b_{I_a} \in \Z_2^n$ so that $I_a \subset supp(s_a)$,
and we call $b_{I_a}s_a$ a {\it signed cycle} of $b_Is$ with the
sign $(-1)^{|I_k|}$. For each signed cycle $b_{I_k}s_k$ with $s_k
= (j_1\cdots j_m)$, the {\it signed cycle-product} of $b_{I_k}s_k$
is the element $g_{j_m}g_{j_{m-1}}\cdots g_{j_1}$ with the sign
$(-1)^{|I_k|}$. For $c\in \G_*$ let $m_i^+(c)$ (resp. $m_i^-(c)$)
be the number of cycles of $x$ whose signed cycle-product lies in
$c$ and has $+$ (resp. $-$) sign. Then $(\rho^+, \rho^-)$ with
$\rho^+(c)=(i^{m_i^+(c)})_{i \geq 1}$ and
$\rho^-(c)=(i^{m_i^-(c)})$ is a pair of partition-valued functions
on $\G_*$ such that $\|\rho^+\|+\|\rho^-\|=n$, and will be called
the {\it type} of the element $(g, b_Is)$. Two elements of $\HG n$
are conjugate if and only if their types are the same. We say that
the conjugacy class $C_{\rho^+, \rho^-}$ is {\it even} (resp. {\it
odd}) if it consists of even (resp. odd) elements. More precisely
if $(g, b_Is)\in C_{\rho^+, \rho^-}$, then $C_{\rho^+, \rho^-}$ is
even (resp. odd) if $|I|$ is even (resp. odd).
\subsection{Split conjugacy classes in $\BG n$}

We can write a general element of $\BG n$ as $(g, z^ka_Is)$ where
\[ z^k a_Is =z^k (a_{I_1}s_1)\cdots (a_{I_q}s_q) ,
\]
and $s=s_1\cdots s_q$ is a cycle decomposition of $s$ and
$I_j\subset supp(s_j)$. We denote by $J^c=\{1,\ldots, n\}-J$ the
complement of a subset $J \subset \{1,\ldots, n\}$.

\begin{lemma} \label{L:conj}
Let $a_Is =(a_{I_1}s_1) \cdots (a_{I_q}s_q)$ be an element of $\B
n$ in its cycle decomposition. Let $J=supp(s_1)\cap I_1^c$, then
\begin{equation*}
(a_Js_1)(a_Is)(a_Js_1)^{-1}=z^{d(s_1)+|I||J|}a_Is.
\end{equation*}
Consequently $(a_{I^c}s)(a_Is)(a_{I^c})^{-1}=z^{d(s)+|I||J|}a_Is$.
\end{lemma}
\begin{proof} Observe that $a_I^2=z^{(|I|+1)|I|/2}$ for any subset
$I$.
For $k\neq 1$ we have
$(a_Js_1)(a_{I_k}s_k)(a_Js_1)^{-1}=z^{|J||I_k|}a_{I_k}s_k$. Therefore
it
reduces to see that
\begin{eqnarray*}
(a_Js_1)(a_{I_1}s_1)(a_Js_1)^{-1}
&=&z^{(|J|+1)|J|/2}a_Ja_{s_1(I_1, J)}s_1 \\
&=&z^{(|J|+1)|J|/2+d(s_1)}a_Ja_{(I_1,J)}s_1\\
&=&z^{(|J|+1)|J|/2+d(s_1)+|J||I_1|}a_J^2a_{I_1}s_1 \\
&=&z^{|J||I_1|+d(s_1)}a_Is_1\, ,
\end{eqnarray*}
where we have used the fact that $supp(s_1)=I_1 \cup J$.
\end{proof}

If two elements of $\BG n$ are conjugate, then clearly their
images are conjugate in $\HG n$. On the other hand, for any
conjugacy class $C$ of $\BG n$, $\theta_n^{-1} (C)$ is either a
conjugacy class of $\BG n$ or it splits into two conjugacy classes
of $\BG n$ (Indeed this holds in a more general setup, cf. e.g.
\cite{Jo1}). A conjugacy class $C_{\rho^+, \rho^-}$ of $\HG n$ is
called {\it split} if the preimage $\theta^{-1}(C_{\rho^+,
\rho^-})$ splits into two conjugacy classes in $\BG n$.
Equivalently, an element $x\in\BG n$ is called {\it split} if $x$
is not conjugate to $zx$ in $\BG n$, then $C_{\rho^+, \rho^-}$ is
split if and only if $\theta^{-1}(C_{\rho^+, \rho^-})$ consists of
split elements. A conjugacy class of $\BG n$ {\it splits} if it
consists of split elements.

The following theorem in the case when $\G=1$ was known in
literature (cf. \cite{Sg, St, Jo2}). It was also stated in
\cite{R} for $\G$ cyclic.

\begin{theorem} \label{T:class}
The conjugacy class $C_{\rho^+, \rho^-}$ in $\HG n$ splits if and only
if

(1) For even $C_{\rho^+, \rho^-}$, we have
$\rho^+\in\mathcal{OP}_n(\G_*)$
and $\rho^-=\emptyset$,

(2) For odd $C_{\rho^+, \rho^-}$, we have $\rho^+=\emptyset$ and
$\rho^-\in \mathcal{SP}^-_n(\G_*)$.
\end{theorem}

\begin{proof}
For $(g, a_Is), (h, a_Jt)\in \BG n$, it follows by definition that

\begin{align} \nonumber
&(h, a_Jt)(g, a_Is)(h, a_Jt)^{-1}\\  \label{E:conj1}
&=(ht(g)(tst^{-1})(h^{-1}), z^{(|J|+1)|J|/2}
a_Ja_{t(I)}a_{tst^{-1}(J)}tst^{-1}).
\end{align}
\noindent($\Longrightarrow$) i) The conjugacy class $C_{\rho^+,
\rho^-}$ is even and split. Suppose on the contrary there is a
part of even integer in $\rho^+(c)$. Without loss of generality we
can assume that $\theta^{-1}(C_{\rho^+, \rho^-})$ contains a
representative element $(g, a_{I}s)$ with the signed cycle
decomposition

$$a_{I}s=(12\cdots r)(a_{I_2}s_2)\cdots (a_{I_p}s_p), $$
where $r =2k$ is even and $I_1$ is empty (we can take all $I_i$
empty corresponding to parts in $\rho^+$) and $|I|$ is even.
Consider the element $$(h, a_Jt)=(h, a_{12\cdots r}(12\cdots
r))\in \BG n, $$ where $h=(h_1, \cdots, h_n)$ with $h_{j}=g_{j}$,
for $j=1, \cdots, r$ and $h_j=1$ otherwise. We claim that

\begin{equation*}
(h, a_Jt)(g, a_Is)(h, a_Jt)^{-1} =z(g, a_Is) .
\end{equation*}
In fact the $j$th component of $ht(g)s(h^{-1})$ equals
$g_{j}g_{j-1}g_{j-1}^{-1}=g_{j}$ for $1\leq j\leq r$ and it also
equals $1\cdot g_j\cdot 1=g_j$ for $j\neq 1, \cdots, r$. Noting
that $st=ts$ we have
\[
ht(g)(tst^{-1})(h^{-1}) = ht(g)s(h^{-1}) =g.
\]
Moreover by Lemma~\ref{L:conj} we have (recall $r=2k$)
\begin{equation*}
(a_{J}t)(a_{I}s)(a_{J}t)^{-1}=z^{(2k-1)+2k\cdot |I|}a_{I}s
=za_{I}s.
\end{equation*}
Thus $(g, a_Is)$ is conjugate to $z(g, a_Is)$ in view of Eqn.
(\ref{E:conj1}). Therefore if the even-parity conjugacy class
$C_{\rho^+, \rho^-}$ splits then $\rho^+\in\mathcal{OP}(\G_*)$.

Now suppose that $\rho^-\neq \emptyset$. Then $\rho^-$ contains at
least two parts since we assume that $C_{\rho^+, \rho^-}$ is even.
Without loss of generality we can assume that
$\theta^{-1}(C_{\rho^+, \rho^-})$ contains an element $(g, a_Is)$
such that
\[
a_{I}s=(a_{i_1}s_1)(a_{i_2}s_2)(a_{I_3}s_3)\cdots (a_{I_p}s_p),
\]
where $i_1\in supp(s_1), i_2\in supp(s_2)$. If both $s_1$ and
$s_2$ are of cycle length $1$, then $(i_1i_2)(g,
a_Is)(i_1i_2)^{-1}=z(g, a_Is)$. Assume that $ord(s_1)\geq 2$, so
$s_1^{-1}(i_1)=i_1'\neq i_1$.
%
    %
%
Consider $(h, a_{i_1}s_1)$, where $h_j=g_{j}$ for $j\in supp(s_1)$
and $h_{j}=1$ otherwise. Then

\begin{eqnarray*}
&& (a_{i_1}s_1)^{-1} a_{I}s (a_{i_1}s_1) \\
 &=& (a_{i_1}s_1)^{-1} (a_{i_1}s_1 a_{i_2}s_2)
   (a_{i_1}s_1) (a_{I_3}s_3)\cdots (a_{I_p}s_p) \\
 &=&  a_{i_2}s_2
   (a_{i_1}s_1) (a_{I_3}s_3)\cdots (a_{I_p}s_p) \\
 &=& z (a_{i_1}s_1) (a_{i_2}s_2) (a_{I_3}s_3)\cdots (a_{I_p}s_p) \\
 &=& za_{I}s.
\end{eqnarray*}
Subsequently $(h, a_{i_1}s_1)^{-1} (g, a_Is) (h, a_{i_1}s_1) =z(g,
a_Is)$. Hence if $C_{\rho^+, \rho^-}$ of even parity splits then
$\rho^-$ is empty. Together with the above we have shown that
split conjugacy class of even parity should have property (1).

ii) The conjugacy class $C_{\rho^+, \rho^-}$ is odd and split. If
on the contrary $\rho^+\neq\emptyset$, we can assume that
$\theta^{-1}(C_{\rho^+, \rho^-})$ contains an element $(g, a_Is)$
with the signed cycle decomposition
\[
a_{I}s=(s_1)(a_{I_2}s_2)\cdots (a_{I_q}s_q),
\]
where $I_1$ is empty and $|I|$ is odd. Take the element $(h,
a_Js_1)$ where $J=supp(s_1)$ and $h_j=g_j$ for $j\in J$ and
$h_j=1$ otherwise. Similarly we can verify that $(h, a_Js_1)(g,
a_Is) (h, a_Js_1)^{-1}=z(g, a_Is)$ by using Lemma~\ref{L:conj}. In
fact
\begin{equation*}
(a_Js_1)(a_Is)(a_Js_1)^{-1}=z^{(|J|-1)+|I||J|}a_Is=za_I,
\end{equation*}
since $|I|$ is odd. Hence $\theta^{-1}(C_{\rho^+, \rho^-})$
does not split if $C_{\rho^+, \rho^-}$ is odd and
$\rho^+\neq\emptyset$.

Next we assume on the contrary that $\rho^-$ contains two
identical parts, then by taking conjugation if necessary we can
assume that $\theta^{-1}(C_{\rho^+, \rho^-})$ contains an element
$(g, a_Is)$ such that
\[
a_{I}s=(a_{i_1}(i_1i_2\cdots i_k))(a_{j_1}(j_1j_2\cdots j_k))\cdots
(a_{I_q}s_q)
\]
and $(g_{i1}, \cdots, g_{i_k})=(g_{i1}, \cdots, g_{i_k})=
(x, 1, \cdots, 1)$ for some $x\in c\in\G_*$.
Consider the element $(1, t)$, where $t=(i_1j_1)\cdots (i_kj_k)$. Then
we
have
\begin{align*}
&(1, t)(g, a_Is)(1, t)^{-1} =(t(g), a_{t(I)}s)\\ &=(g,
a_{j1}(j_1\cdots j_k)a_{i_1}(i_1\cdots i_k)\cdots)\\ &=(g,
za_{i_1}(i_1\cdots i_k)a_{j_1}(j_1\cdots j_k)\cdots)=(g, za_Is),
\end{align*}
which is a contradiction.

\noindent($\Longleftarrow$) Suppose that (1) holds. If on the
contrary the even conjugacy class $C_{\rho^+, \emptyset}$ does not
split, then we can assume that $\theta^{-1}(C_{\rho^+,
\emptyset})$ contains an element $(g, s)$, where $s=s_1\cdots s_q$
with each $s_i$ being odd cycle, and $(h, a_Jt)(g, s)(h,
a_Jt)^{-1}=z(g, s)$. Therefore $(a_Jt)s(a_Jt)^{-1}=zs$, i.e.,
$za_{s(J)}=a_J$, which in particular implies that $supp(s)\subset
J$. Then $a_{s(J)}=z^{d(s)}a_J$ by Eq.~(\ref{eq_perm}), and so
$d(s) =1 \bmod{2}$ which contradicts with the assumption on $s$.

Suppose that (2) holds. Assume on the contrary that the odd
conjugacy class $C_{\emptyset, \rho^-}$ does not split. Since an
identification of two elements in $\HG n$ implies that their
respective components in $\HG{|\rho^-(c)|}$ are already equal, we
can assume that $\rho^-$ consists of one strict partition
$\rho(c)^-$ for some $c$. Thus $\theta^{-1}(C_{\emptyset,
\rho^-})$ contains a non-split element $(g, a_Is)=(g,
a_{i_1}s_1\cdots a_{i_q}s_q)$, where $q$ is odd and $i_k\in
supp(s_k)$. Let $(h, a_Jt)(g, a_Is)(h, a_Jt)^{-1}=z(g, a_Is)$ for
some element $(h, a_Jt)$. It follows that $t$ commutes with $s$
and $s$ is a product of disjoint cycles with mutually distinct
orders, the permutation $t$ equals $s_1^{r_1}s_2^{r_2}\cdots
s_q^{r_q}$ for $0\leq r_i\leq ord(s_i)$. Thus we can write $a_Jt
=(a_{J_1}t_1)\cdots (a_{J_q}t_q)$ with $t_k=s_k^{r_k}$. As in the
proof of Lemma \ref{L:conj} we have

\begin{align*}
&(a_Jt)(a_{i_1}s_1)(a_Jt)^{-1} \\
 &=(a_{J_1}t_1) (a_{J_2}t_2) \cdots (a_{J_q}t_q)
 (a_{i_1}s_1) ((a_{J_2}t_2) \cdots (a_{J_q}t_q))^{-1}
(a_{J_1}t_1)^{-1}\\
&=z^{|J|-|J_1|}(a_{J_1}t_1)(a_{i_1}s_1)(a_{J_1}t_1)^{-1},
\end{align*}
which must equal $a_{i_1}s_1$ up to a power of $z$. Set
$(a_{J_1}t_1)(a_{i_1}s_1)(a_{J_1}t_1)^{-1} =z^*a_{i_1}s_1$ where
$*$ is $0$ or $1$. We claim that $*$ is always $0$. Without loss
of generality we let $J_1 =\{1, 2,\cdots, r\}$, $i_1=1$,
$s_1=(12\cdots k)$ with $0\leq r\leq k$, then
\begin{equation*}
a_{J_1}a_{s_1^{r_1}(1)}a_{s_1(J_1)}z^{(r+1)r/2}
=a_{12\cdots r}a_{r_1}a_{23\cdots r+1}z^{(r+1)r/2}= z^{*}a_1
\end{equation*}
implies that $r_1=r+1$, which in turn implies the exponent $*$ is
equal to $0$. Therefore
$(a_Jt)(a_{i_1}s_1)(a_Jt)^{-1}=z^{|J|-|J_1|}a_{i_1}s_1$ and
similarly we have
\begin{align*}
&(a_Jt)(a_{I}s)(a_Jt)^{-1}=(a_Jt)(a_{I_1}s_1)\cdots
(a_{I_q}s_q)(a_{J}t)^{-1}\\ &=z^{q|J|-(|J_1|+\cdots
|J_q|)}a_{I}s=z^{(q-1)|J|}a_{I}s =a_{I}s,
\end{align*}
since $q$ is odd. This is a contradiction.
\end{proof}

For $\rho\in\mathcal{OP}_n(\G_*)$ we let $D_{\rho}=D_{\rho}^+$ be
the split conjugacy class in $\BG n$ containing the elements
$(g,s)$ of type $(\rho, \emptyset)$ with $g \in \G^n, s \in S_n$.
Then $zD_{\rho}^+$ is the other conjugacy class lying in
$\theta_n^{-1}(C_{\rho, \emptyset})$, which will be denoted by
$D_{\rho}^-$.

For each partition-valued function $\rho=(\rho(c))_{c\in\G_*}$ we
define $$ Z_{\rho}=\prod_{c\in\G_*, i\geq 1}m_i(c)!
i^{m_i(c)}\zeta_c^{l(\rho(c))}, $$ which is the order of the
centralizer of an element of conjugacy type $\rho=( \rho(c))_{ c
\in \G_*}$ in $\Gn$ (cf. \cite{M}). It follows from Theorem
\ref{T:class} that the order of the centralizer of an element in
the conjugacy class $D_{\rho}^+$ is
\begin{equation} \label{E:centralizer}
2^{1+l(\rho)}Z_{\rho}.
\end{equation}
\section{Spin supermodules over $\BG n$} \label{sec_spinsuper}
In this section, we first recall some general facts about {\em
spin} supermodules over a superalgebra. This will be applied to
(the group algebra of) $\BG n$. We introduce the Grothendieck
group $R^-(\BG n)$ of spin supermodules over $\BG n$, and
construct the so-called basic spin supermodules of $\BG n$. We
further define a natural Hopf algebra structure on $\tRG =
\oplus_{n \geq 0} R^-(\BG n)$, and introduce a weighted bilinear
form on $\tRG$ associated to any given self-dual virtual character
$\wt$ of $\G$. When $\wt$ is trivial, the weighted bilinear form
reduces to the standard one.
\subsection{Superalgebras and supermodules}
A complex {\it superalgebra} $A=A_0\oplus A_1$ is a $\mathbb
Z_2$-graded complex vector space with a binary product $A\times
A\longrightarrow A$ such that $A_iA_j\subset A_{i+j}$. A $\mathbb
Z_2$-graded vector space $V=V_0\oplus V_1$ is a {\it supermodule}
for a superalgebra $A$ if $A_iV_j\subset V_{i+j}$. A linear map
$f: M\rightarrow N$ between two $A$-supermodules is a morphism of
degree $i$ if $f(M_j)\subset M_{i+j}$ and for any homogeneous
element $a\in A$ and any homogeneous vector $m\in M$ we have
\[ f(am)=(-1)^{p(a)p(f)}af(m).
\]
Let $V=V_0\oplus V_1$ and $W=W_0\oplus W_1$ be two supermodules.
The tensor product $V\otimes W$ is also a supermodule with
$(V\otimes W)_i=\sum_{k+l=i (mod 2)}V_k\otimes W_l$. The notions
of submodules, irreducible supermodules etc are defined similarly
as usual.

Let $\mathbb C^{r|s}$ be the $\Z_2$-graded vector space $\mathbb
C^r\oplus \mathbb C^s$. The algebra $M(r|s)$ consisting of all
linear transformations on $\mathbb C^{r|s}$ inherits a natural
$\Z_2$ grading from $\mathbb C^{r|s}$. It is easily seen that
$M(r|s)$ is a simple superalgebra.

Another example of simple superalgebra is $Q(n)$, which is the
subalgebra of $M(n|n)$ consisting of matrices of the form

\begin{align*}
\begin{bmatrix} C & D\\ D & C\end{bmatrix}, \quad C ,D \in
M_n(\mathbb C).
\end{align*}
The left multiplication of $Q(n)$ on $\C^{n|n}$ gives an
irreducible supermodule structure on it.

A well-known result due to C.~T.~C.~Wall says that these
superalgebras are the only simple superalgebras over $\mathbb C$.
In the sequel we will say the supermodules $\mathbb C^{r|s}$ and
$\C^{n|n}$ are of type $M$ and $Q$ respectively.

Let $G$ be a finite group and let $p: G\rightarrow \Z_2$ be a
group epimorphism. We denote by $G_0$ the kernel of $p$ which is a
subgroup of $G$ of index $2$. We regard $p(\cdot)$ as a parity
function on $G$ by letting the degree of elements in $G_0$ be $0$
and letting the degree of elements in the complementary $G_1 = G
\backslash G_0$ be $1$. Elements in $G_0$ (resp. $G_1$) will be
called even (resp. odd). In addition we assume that $G$ contains a
distinguished even central element $z$ of order $2$. A {\em spin}
supermodule over $G$ is a supermodule over the group superalgebra
$\C[G]$ such that $z$ acts as $-1$. The group superalgebra is
semisimple (cf. \cite{Jo1}), i.e. decomposes into a direct sum of
simple superalgebras. We will refer the corresponding supermodules
as of type $M$ and type $Q$.

Now let us return to our main example $\BG n$. It is easy to see
that the characters of spin supermodules vanish on nonsplit
classes. Let $(-1)^p$ be the one-dimensional representation of
$\BG n$ given by $x\longmapsto (-1)^{p(x)}$. A representation
$\pi$ of $\BG n$ is called a {\it double spin} representation if $
(-1)^{p}\pi \simeq\pi$. If $\pi'=(-1)^{p}\pi \ne\pi$, then $\pi'$
and $\pi$ are called {\it associate spin} representations of $\BG
n$. By the general theory of supermodules (cf. \cite{Jo1}) and
Theorem~\ref{T:class} we obtain the following proposition.

\begin{proposition}
The number of irreducible double spin representations over $\BG n$
is equal to $|{\mathcal SP}_n^+(\G_*)|$, and the number of pairs
of irreducible associate spin representations is $|{\mathcal
SP}_n^-(\G_*)|$. The number of irreducible spin supermodules of
$\Gn$ is $|\mathcal{SP}_n(\G_*)|$.

If $V$ is an irreducible $\BG n$-supermodule of type $M$, then its
underlying module $|V|$ (by forgetting the $\Z_2$-grading
structure) is an irreducible double spin $\BG n$-module. If $P$ is
an irreducible $\BG n$-supermodule of type $Q$, then $|P|\simeq
N\oplus N'$ where $N$ and $N'$ are a pair of irreducible associate
spin $\BG n$-modules.
\end{proposition}

\begin{remark}
It has been observed (cf. \cite{St, Jo2, Naz}) that there is a
very close connection between the representations of $\B n$ and
those of a spin cover $\tS n$ of the symmetric group $S_n$.
Yamaguchi \cite{Y} explains such an phenomenon in an elegant way
by establishing an isomorphism between the group superalgebra $\C
[\B n]/ \langle z=-1 \rangle$ and the (outer) tensor product of
the group superalgebra $\C[\tS n]/ \langle z=-1 \rangle$ with the
complex Clifford algebra of $n$ variables. (Note that a Clifford
algebra admits a unique irreducible supermodule). In view of this,
there is also an isomorphism by substituting $\B n$ and $ \tS n$
with $\BG n$ and $\tG n$ respectively. This isomorphism provides a
direct connection between the constructions in the present paper
and \cite{FJW2}.
\end{remark}

Given a $\Z_2$-graded finite group $G$ and a $\Z_2$-graded
subgroup $H$ that contain an even central element $z$ of order
$2$, we can define the induction and restriction of supermodules
similarly as usual. In particular the induced supermodule (or a
restriction) of a spin supermodule remains to be a spin
supermodule. The Mackey theorem remains true in this setup (cf.
\cite{FJW2}, Sect. 3).
\subsection{The Grothendieck group $R^-(\BG n)$}
Let $R(\G)=\oplus_{i=0}^{r} \mathbb C\g_i$ be the space of
complex-valued class functions on $\G$. The usual bilinear form on
$R(\G )$ is defined as follows:
\begin{eqnarray*}
\langle f, g \rangle_{\G} = \frac1{ | \G |}\sum_{x \in \Gamma}
          f(x) g(x^{ -1})
 = \sum_{c \in \Gamma_*} \zeta_c^{ -1} f(c) g(c^{ -1}).
\end{eqnarray*}
Then $\Rz :=\oplus_{i=0}^{r} \mathbb Z\g_i$ becomes an integral
lattice in $R( \G )$ under the standard bilinear form: $\langle
\g_i, \g_j \rangle =\delta_{ij}$.

A {\it spin class function} on $\BG n$ is class function $f: \BG
n\longrightarrow \mathbb C$ such that $f(zx)=-f(x)$, hence it
vanishes on non-split conjugacy classes. A {\it spin super class
function} $f$ on $\BG n$ is a spin class function such that it
vanishes further on odd conjugacy classes. In other words, $f$
corresponds to a complex functional on $\mathcal{OP}_n(\G_*)$ in
view of of Theorem~\ref{T:class}.

Let $R^-(\BG n)$ be the $\mathbb C$-span of spin super class
functions on $\BG n$. The standard bilinear form on $R^-(\BG n)$
is given by
\begin{equation*}
\langle f, g \rangle_{\BG n}=\frac1{ | \BG n |}\sum_{x \in \BG n}
          f(x) g(x^{ -1})
= \sum_{\rho\in \mathcal{OP}_n(\G_*)} \frac
1{2^{l(\rho)}Z_{\rho}}f(\rho)
g(\overline{\rho}),
\end{equation*}
where $f, g\in R^-(\BG n)$, $f(\rho)=f(D_{\rho}^+)$, and Eqn.
(\ref{E:centralizer}) is used here. The following holds for a
general $\Z_2$-graded finite group \cite{Jo1}.


\begin{proposition} \label{P:criterion}
The characters of irreducible spin supermodules over $\BG n$ form
a $\mathbb Z$-basis of $R^-(\BG n)$. Let $\phi$ and $\ga$ be two
irreducible characters of spin supermodules, then

\begin{equation}\label{E:superorth}
\langle \phi, \ga\rangle=\begin{cases} 1 & \mbox{if $\phi\simeq
\ga$, type $M$}\\ 2 & \mbox{if $\phi\simeq \g$, type $Q$}\\ 0 &
\mbox{otherwise}\end{cases}.
\end{equation}
Conversely, if $\langle f, f\rangle=1$ for $f\in R^-(\BG n)$, then
$f$ or $-f$ affords an irreducible spin $\BG n$-supermodule of
type $M$.
\end{proposition}
\subsection{Basic spin supermodules}
Let $L_n$ be the Clifford algebra generated by $e_1, e_2, \cdots
e_{n}$ with relations:
\begin{equation}\label{E:clifford}
\{e_i, e_j\} =e_ie_j+e_je_i=-2\delta_{ij}.
\end{equation}
The group $\B n$ acts on $L_n$ by $a_i.e_I=e_ie_{I}$ and
$s.e_i=e_{s(i)}$.

The Clifford algebra $L_n$ has a natural $\Z_2$-grading given by
the parity $p$ such that $p(e_i)=1$. The set $\{a_I\}, I\subset
\{1, \ldots, n\}$ form a basis of $L_n$. We observe that
$s(a_I)=(-1)^{d(s)}a_I$ $(s \in S_n)$ if and only if $I$ is a
union of $supp(s_i)$, where $s=s_1\cdots s_l$ is a cycle
decomposition.

\begin{proposition}\cite{Jo2} The module $L_n$ is an irreduicble $\B
n$-supermodule. The value of its character $\chi_n$ at the
conjugacy class $D_{\alpha}^+$ is given by
\begin{equation}\label{E:spinchar}
\chi_n(\alpha)=2^{l(\alpha)}, \qquad
\mbox{$\alpha\in\mathcal{OP}_n$}.
\end{equation}
\end{proposition}

Let $V$ be a $\G$-module afforded by the character $\g \in R(\G)$,
the tensor product $V^{\otimes n}$ is a $\Gn$-module by the direct
product action of $\G^n$ combined with permutation action of the
symmetric group $S_n$. More explicitly the action is given by
\begin{equation*}
(g, s).(v_i\otimes\cdots\otimes v_n) =g_1v_{s^{-1}(1)}\otimes
\cdots g_nv_{s^{-1}(n)},
\end{equation*}
where $g=(g_1, \cdots, g_n)\in \G^n, s\in S_n$. We denote the
resulting character by $\eta_n(\ga)$. Indeed one can extend (cf.
\cite{W, FJW1}) this construction to define a map $\eta_n: R (\G)
\rightarrow R (\Gn)$.

Let $V$ be a $\G$-module afforded by the character $\g \in R(\G)$,
and let $U$ be a spin supermodule of $\B n$ with the character
$\pi_n$. The tensor product $V^{ \otimes n}\otimes U$ has a
canonical spin supermodule structure for $\BG n$ as follows. For
any $g=(g_1, \cdots, g_n)\in \G^n$, and an element $(g, a_Is)$ in
$\BG n$, the action is given by
\begin{equation*}
(g, a_Is).(v_i\otimes\cdots\otimes v_n\otimes u)
=g_1v_{s^{-1}(1)}\otimes \cdots g_nv_{s^{-1}(n)} \otimes (a_Isu).
\end{equation*}
It is easy to see that if $V$ and $U$ are irreducible, then so is
$V^{ \otimes n}\otimes U$. We denote by $\pi_n(\g)$ the character
of this spin supermodule. The following result can be proved as in
\cite{FJW1, FJW2} for similar results.

\begin{proposition}
Let $\pi_n$ be the character of a spin $\B n$-supermodule. Then
the value of the character ${\pi_n}(\g)$ at an element $(g,
a_Is)\in D_{\rho}^+$ in the even split conjugacy class
$D_{\rho}^{+}$ is given by

\begin{equation} \label{eq_term}
\pi_n(\g)((g, a_Is)) =\pi_n(a_Is)\prod_{c\in\G_*}
\g(c)^{l(\rho(c))}.
\end{equation}
\end{proposition}

The map $\pi_n$ can be extended $R(\G)$ to $R^-(\BG n)$ (compare
\cite{FJW2}):
\begin{equation}  \label{eq_virt}
  \pi_n (\beta - \g) =
  \sum_{m =0}^n ( -1)^m Ind_{\BG {n -m} \hat{\times} \BG m }^{\BG n}
   [ \pi_{n -m} (\beta) \otimes \pi_m (\g ) ] .
\end{equation}
In particular, when $U$ is the spin $\B n$-supermodule $L_n$ we
denote by $\chi_n(\g)$ the character of the $\BG n$-supermodule
$L_n(V) := V^{\otimes n}\otimes L_n$. We will refer to $L_n(V)$
(associated to irreducible $V$) as the {\em basic} spin
supermodules over $\BG n$.

\begin{corollary} \label{C:char}
The character values of $\chi_n(\g)$ on the conjugacy classes
$D_{\rho}^{\pm}$  are given by
\begin{equation}\label{E:charvalue}
\chi_n(\g)(D_{\rho}^{\pm})= \pm 2^{l(\rho)} \prod_{c\in\G_*}
\g(c)^{l(\rho(c))}, \quad\rho\in \mathcal{OP}_n(\G_*).
\end{equation}
\end{corollary}
\subsection{Hopf algebra structure on $\tRG$}\label{S:Hopf}
Our main object of this paper is the space
\[
  \tRG = \bigoplus_{n\geq 0} R^-(\BG n).
\]

Let $\BG n\tilde{\times}\BG m$ be the direct product
of $\BG n$ and $\BG m$ with a twisted multiplication
$(t, t')\cdot (s, s')=(ts, z^{p(t')p(s)}t's'),
$
where $s, t\in\BG n, s', t'\in\BG m$
are homogeneous. We define the {\it spin product}
of $\BG n$ and $\BG m$ by
\begin{equation}
\BG n\hat{\times}\BG m={\BG n\tilde{\times}\BG m}/
\{(1, 1), (z, z)\},
\end{equation}
which can be embedded into the spin group $\BG{n+m}$ canonically
by letting $ (t_i', 1)\mapsto t_i, (1, t_j'')\mapsto t_{n+j},$
$i=1, \ldots, n-1$, $j=1,\ldots, m-1$. We will identify $\BG
n\hat{\times}\BG m$ with its image in $\BG {m+n}$. The subgroup
$\BG n\hat{\times}\BG m$ has a distinguished subgroup of index $2$
consisting of even elements with $p=0$. Therefore we can define
$R^-(\BG n\hat{\times}\BG m)$ to be the space of super spin class
functions.

For two spin modules $U$ and $V$ of $\BG n$ and $\BG m$ we define
the {\it (outer)-tensor product} by
\begin{equation*}
(t, s)\cdot (u\otimes v)=(-1)^{p(s)p(u)}(tu\otimes sv),
\end{equation*}
where $s$ and $u$ are homogeneous elements. Since $(z, 1)\cdot
(u\otimes v)=-u\otimes v$, the tensor $U\otimes V$ is a spin $\BG
n\hat{\times}\BG m$-supermodule. The following is a
straightforward generalization of a result in \cite{Jo2} for $\B
n$.

\begin{proposition}\label{P:supermult}
Let $U$ and $V$ be simple
supermodules for $\BG n$ and $\BG m$ respectively. Then
\item 1) If both $U$ and $V$ are of type M, then $U\otimes V$
is a simple $\BG n\hat\times\BG m$-supermodule of type M.
\item 2) If $U$ and $V$ are of different type, then
$U\otimes V$
is a simple $\BG n\hat\times\BG m$-supermodule of type Q.
\item 3) If both $U$ and $V$ are type Q, then $U\otimes V
\simeq N\oplus N$ for some simple $\BG n\hat\times\BG
m$-supermodule $N$ of type M. We will denote $N$ by $2^{-1}U
\otimes V.$
\end{proposition}

For a simple supermodule $V$ we define
$c(V)=0$ if $V$ is type M and $c(V)=1$ if $V$ is type Q.
Set $c(V_1, V_2)=c(V_1)c(V_2)$ for two simple
supermodules of $\mathbb C[\BG n]$ and $\mathbb C[\BG m]$
respectively.
It is easy to see that $c$ satisfies the cocycle condition
\begin{equation} \label{E:cocycle}
c(V_1, V_2)+c(V_1\otimes V_2, V_3)
=c(V_2, V_3)+ c(V_1, V_2\otimes V_3).
\end{equation}
It follows from Proposition \ref{P:supermult} and Eqn.
(\ref{E:cocycle}) that the tensor product defines an isomorphism:
\begin{equation*}\label{E:twprod}
R^-(\BG n ) \bigotimes R^-(\BG m)
 \stackrel{\phi_{n, m}}\longrightarrow R^-(\BG n \hat{\times} \BG m),
\end{equation*}
where $\phi_{n, m}(V_1\otimes V_2)=2^{-c(V_1, V_2)}V_1\otimes V_2$ for
simple
modules $V_1$ and $V_2$.

We now define a multiplication on $\tRG$ by the composition

\begin{equation}\label{E:hopf1}
 m: R^-(\BG n ) \bigotimes R^-(\BG m)
 \stackrel{\phi_{n, m}}\longrightarrow R^-(\BG n \hat{\times} \BG m)
 \stackrel{Ind}{\longrightarrow} R^-( \BG{n + m}),
\end{equation}
and a comultiplication on $\tRG$ by the composition

\begin{eqnarray}\label{E:hopf2}
\Delta: R^-(\BG n ) &\stackrel{Res }{\longrightarrow}
 &\bigoplus_{m =0}^n R^-( \BG {n - m}\times \BG m) \\
 &\stackrel{\phi^{-1}}{\longrightarrow}&
 \bigoplus_{m =0}^n R^-( \BG {n - m}) \bigotimes
R^-(\BG m). \nonumber
\end{eqnarray}
Here $Ind$ and $Res$ denote the induction and restriction of spin
supermodules respectively, and the isomorphism $\phi^{-1}$ is
given by $\oplus_{m=0}^{n}\phi^{-1}_{n-m, m}$.

\begin{theorem}
The above operations define a Hopf algebra structure for $\tRG$.
\end{theorem}
\begin{proof} The associativity follows from
the obervation that
the two different embeddings are conjugate:
\[  (\BG n\hat{\times}\BG m)\hat{\times}\BG l\hookrightarrow
\BG{n+m+l}\hookleftarrow \BG n\hat{\times}
(\BG m\hat{\times}\BG l).
\]
Using the cocycle condition (\ref{E:cocycle}) we can check the
coassociativity as we did in \cite{FJW2}. Using the cocycle $c$
again and super analog of Mackey's theorem we can check that
$\Delta$ preserves the multiplication structure, for details see
\cite{FJW2} in a similar situation (also compare \cite{Z, W} for
Hopf algebra structures in different but related setups).
\end{proof}
\subsection{A weighted bilinear form on $R(\G)$ and $R^-(\BG n)$}
Let $\wt$ be a self-dual virtual character in $R_{\Z}(\G)$, i.e.
$\wt(c)=\wt(c^{-1})$, the {\em weighted} bilinear form on
$R_{\Z}(\G)$ (cf. \cite{FJW1}) is defined by $$
  \langle f, g \rangle_{\wt } = \langle \wt * f ,  g \rangle_{\G },
   \quad f, g \in R_{\Z}( \G),
$$ where $*$ the product of two characters. The self-duality of
$\xi$ is equivalent to the condition that the matrix of the
weighted bilinear form is symmetric. When $\xi$ is the trivial
character $g_0$ of $\G$, the weighted bilinear form becomes the
standard one. One extends the weighted bilinear form to $R(\G)$ by
bilinearity.

Let $V =V_0 +V_1$ be a spin supermodule for $\BG n$ and $W$ a
module for $\Gn$. The tensor product $W\circledast V=W\otimes
V_0\oplus W\otimes V_1$ carries a natural $\Z_2$-grading and
admits a natural spin $\BG n$-supermodule structure by letting
\begin{equation}\label{E:circtensor}
(g, a_I \sigma )(w\otimes v)=(g, {\sigma}) \cdot w\otimes (g, a_I
\sigma)\cdot v,
\end{equation}
where $g\in \G^n, a_I \in \Pi_n, \sigma\in S_n.$ This gives rise
to a morphism:
\begin{equation}\label{E:tensorprod}
R(\Gn)\otimes R^-(\BG n)\stackrel{\circledast}{\longrightarrow}
R^-(\BG n).
\end{equation}

Recall we have defined $\eta_n(\xi) \in R(\Gn)$ associated to $\xi
\in R(\G)$. Its character value at the class $\rho=(\rho(c))_{c
\in \G_*}$ is given by (cf. \cite{FJW1, M})
\begin{equation}\label{E:weightchar}
\eta_n(\xi)(\rho)=\prod_{c\in\G_*}\xi(c)^{l(\rho(c))}.
\end{equation}
Thus $\eta_n(\xi)$ is self-dual as long as $\xi$ is.

We now introduce a {\em weighted bilinear form} on $R^-( \BG n)$ by
letting
\begin{align}\nonumber
  \langle  f, g\rangle_{\wt} &=
   \langle \eta_n (\wt ) \circledast f, g \rangle_{\BG n}
     \\ \label{E:innerprod}
&=\sum_{\rho\in\mathcal{OP}_n(\G^*)}\frac1{2^{l(\rho)}Z_{\rho}}f(\rho)
g(\overline{\rho})\prod_{c\in\G_*}\xi(c)^{l(\rho(c))},
\end{align}
where $f, g \in R^-( \BG n)$. The self-duality of $\eta_n (\wt)$
implies that the bilinear form $\langle \ , \ \rangle_{\wt}$ is
symmetric. When $\xi$ is taken to be the trivial character $\g_0$,
then it reduces to the standard bilinear form on $R^-( \BG n)$.

The bilinear form on $\tRG = \bigoplus_{n} R^-(\BG n)$ is given by
\[
\langle u, v \rangle_{\wt}
 = \sum_{ n \geq 0} \langle u_n, v_n \rangle_{\wt, \BG n } ,
\]
where
$u = \sum_n u_n$ and $v = \sum_n v_n$ with $u_n, v_n\in R^-(\BG n)$.

\begin{remark}
When $\G$ is finite subgroup of $SL_2(\mathbb C)$, an importance
choice for $\xi$ is $\wt = 2 \g_0 - \pi,$ where $\pi$ is the
character afforded by the $2$-dimensional representation of $\G$
given by the embedding of $\G$ in $SL_2(\mathbb C)$. We shall see
that the weighted bilinear form $\langle \ , \ \rangle_{\wt}$ on
$\tRG$ becomes positive semi-definite. This will play an important
role in the later part of this paper.
\end{remark}
\section{Identification of $\tRG$ as a Fock space}\label{sec_heis}
In this section, we first recall a twisted Heisenberg algebra
$\thg$ and its Fock space $\tSG$ together with a bilinear form. We
define an action of $\thg$ on $\tRG$ in terms of group-theoretic
maps. We further show that there is a natural isometric
isomorphism from $\tRG$ to $\tSG$ which is compatible with the
Hopf algebra structure and Heisenberg algebra action on both
spaces.
\subsection{A twisted Heisenberg algebra $\thg$}
Associated with a finite group $\Gamma$ and a self-dual class
function $\wt\in R(\G)$, a twisted Heisenberg algebra $\thg $ (cf.
\cite{FJW2}) is generated by $a_m(\gamma), m \in 2\mathbb Z+1,
\gamma \in\G^*$ and a central element $C$, subject to the
relations:
\begin{equation}  \label{eq_heis}
[a_m( \gamma), a_n(\gamma ')]
 = \frac m2 \delta_{m, -n} \langle \gamma, \gamma' \rangle_{\wt } C,
 \quad m, n \in 2\mathbb Z+1, \, \g, \g ' \in \G^*.
\end{equation}
For
$\g = \sum_{ i =0}^r s_i \g_i \in R(\G )$ $(s_i \in \mathbb C)$
we write
$ a_m ( \g ) = \sum_i s_i \, a_m (\g_i )$.
The center of $\thg $ is spanned by $C$ together with
$a_m ( \g),  m \in 2\mathbb Z+1, \g \in R_0$,
the radical of the bilinear form
$\langle \cdot , \cdot \rangle_{\wt }$ in $R(\G)$.

We introduce another basis for
$\thg$:
\begin{equation}\label{E:classba}
 a_{ m}( c) = \sum_{ \g\in \G^*} \gamma(c^{-1}) a_m( \g ),
\quad m\in 2\mathbb Z+1, c \in \G_*.
\end{equation}
Equivalently $a_m( \g )
   = \sum_{c \in \G_*} \zeta_c^{ -1}
   \g  (c) a_m(c).$
Then the commutation relations (\ref{eq_heis}) imply that
 \begin{eqnarray} \label{prop_orth}
  [ a_m( {c'}^{ -1}), a_n( c )]
    & =& \frac m2 \delta_{m, -n}\delta_{c', c} \zeta_c \wt (c) C,
   \quad c, c' \in \G_*,
\end{eqnarray}
where $m, n\in Z$.
\subsection{The Fock space $\tSG$}
The Fock space $\tSG $ is defined to be the symmetric algebra
generated by $a_{-n}(\g), n \in 2\mathbb Z_++1, \g\in \Gamma^*$.
There is a natural grading on $\tSG $ by letting
 $$ \deg (a_{ -n}( \g)) = n , \qquad n\in 2\mathbb Z_++1,
$$
which makes $\tSG $ into a $\mathbb Z_+$-graded space.

We define an action of $\thg$ on $\tSG$ as follows: $a_{-n}( \g),
n>0$ acts as multiplication operator on $ \tSG $ and $C$ as the
identity operator; $a_n (\g),$ $ n >0$ acts as a derivation of
algebra:
\begin{align*}
  a_n (\g).& a_{-n_1}( \alpha_1)
    \ldots a_{-n_k}( \alpha_k) \\
 &= \sum_{i =1}^k \delta_{n,n_i}
 \langle \g , \alpha_i \rangle_{\wt }
   a_{-n_1}( \alpha_1)\ldots
  \check{a}_{-n_i}( \alpha_i) \ldots a_{-n_k}(\alpha_k )  .
\end{align*}
Here $n_i> 0, \alpha_i \in R(\G)$ for $i =1, \ldots , k$, and
$\check{a}_{-n_i}( \alpha_i)$ means the very term is deleted. Note
that $\tSG$ is not an irreducible representation over $\hg$ in
general since the bilinear form $\langle \ , \ \rangle_{\wt}$ may
be degenerate.

Denote by $\SGO$ the ideal in the symmetric algebra $ \tSG$
generated by $a_{-n}(\g), n \in \mathbb N, \g \in R_0$. Denote by
$\tSGG$ the quotient $\tSG /\SGO$. It follows from the definition
that $\SGO$ is a submodule of $\tSG$ over the Heisenberg algebra
$\thg$. In particular, this induces a Heisenberg algebra action on
$\tSGG$ which is irreducible. The unit $1$ in the symmetric
algebra $\tSG$ is the highest weight vector. We will also denote
by $1$ its image in the quotient $\tSGG$.
\subsection{The bilinear form on $\tSG $}
The Fock space $ \tSG $ admits a bilinear form $\langle \ ,  \
\rangle_{\wt } '$ determined by
\begin{eqnarray}  \label{eq_bili}
 \langle 1, 1 \rangle_{\wt}' = 1, \quad
 a_n(\g)^* = a_{-n}(\g), \qquad n\in 2\mathbb Z+1, n \in \G^*.
\end{eqnarray}
Here $a_n(\g)^*$ denotes the adjoint of $a_n(\g)$.

For $\la\in\mathcal{OP}$ we write $a_{-\la}( \g) = a_{-\la_1}(
\g)a_{ - \la_2}( \g) \cdots$ for $\g \in \G^*$, and $a_{ - \la} (c
)  = a_{ - \la_1}(c) a_{ - \la_2} (c) \cdots$ for $c\in\G_*$. For
$\lambda = ( \lambda (\g) )_{ \g \in \G^*} \in {\mathcal
{OP}}(\G^* )$, we define $$
  a_{ - \lambda} = \prod_{\g \in \G^*}  a_{ - \lambda(\g)}(\g).
$$
and similarly $a_{- \rho}' = \prod_{ c \in \G_*} a_{ - \rho (c)} (c)$
for $\rho = ( \rho (c) )_{ c \in \G_* } \in
 \mathcal {OP} ( \G_* )$.
It is clear that both $\{a_{ - \lambda}\}$ and $\{a_{-\rho}'\}$ are
$\mathbb C$-bases for $\tSG $.

It follows from Eqn. (\ref{prop_orth}) and (\ref{eq_bili}) that
\begin{equation}  \label{eq_inner}
  \langle a_{ - \rho'}', a_{ - \overline{\rho} }' \rangle_{\wt }'
  = \delta_{\rho ', \rho }
    \frac{Z_{\rho}}{2^{l(\rho)}}
\prod_{c \in \G_*} \wt (c)^{l (\rho (c))},
 \quad \rho ', \rho  \in \mathcal{OP}(\G_*).
\end{equation}
The bilinear form $\langle \ , \ \rangle_{\wt} '$
 induces a bilinear form on $\tSGG$ which will be denoted
by the same notation.
\subsection{The characteristic map $\ch$}
We define the {\it characteristic} map $ch: \tRG \longrightarrow
\tSG$ by letting (compare with \cite{FJW2})
\begin{equation}\label{E:ch}
ch (f)
= \sum_{\rho \in \mathcal {OP}(\G_*)}\frac{1}{Z_{\rho}}
f_{\rho}
a_{-\overline{\rho}}',
\end{equation}
where $f_{\rho}=f(D_{\rho}^+)$. In the case when $\G$ is trivial,
this is essentially the same as defined in \cite{Jo2}, and a
different approach is given in \cite{W2}.

For $c\in\G^*$ and $n\in 2\mathbb Z+1$, we let
$c_n (c \in \G_*)$ be the split conjugacy class $D^+_{\rho}$ in $\BG
n$
such that $\rho(c)=(n)$ and $\rho(c')=\emptyset$ for
$c'\neq c$.
Denote by $\sigma_n (c )$ the super class function on $\BG n$ which
takes
value $n\zeta_c$
on elements in the conjugacy class $c_n$
and $0$ elsewhere. For
$\rho = \{(r^{m_r (c)}) \}_{r \geq 1, c \in \G_*}
\in \mathcal {OP}_n (\G_*)$,
$\sigma_{\rho} =
\prod_{r \geq 1, c \in \G_*} \sigma_r (c)^{m_r (c)}$
is the class function on split conjugacy classes in $\BG n$ which
takes value
$Z_{\rho}$ on the conjugacy class $D_{\rho}^+$ and
$0$ elsewhere. Given $\g \in R(\G)$, we
denote by $\sigma_n (\g )$ the class function on $\BG n $ which takes
value $n\g (c) $ on elements in the class $c_n, c \in \G_*$,
and $0$ elsewhere.

The following lemma follows from definitions.
\begin{lemma}  \label{lem_isom}
  The map $ch$ sends $\sigma_{\rho}$ to $a_{ - \rho} '$.
 In particular, it sends $\sigma_n (c)$ to $a_{ -n} (c)$ in $\tSG$
 while sending $\sigma_n(\g )$ to $a_{ -n} ( \g )$.
\end{lemma}
\subsection{Action of $\thg$ on $\tRG$}
We define $ \widetilde{a}_{ -n} (\gamma), n\in 2\mathbb Z_++1 $ to
be a map from $\tRG$ to itself by the following composition
\[  R^- (\BG m) \stackrel{ \sigma_n ( \g ) \otimes }{\longrightarrow}
  R^-(\BG n) \bigotimes R^- (\BG m)  \stackrel{{Ind}
}{\longrightarrow}
  R^- ( \BG{n +m}).
\]
We also define $ \widetilde{a}_{ n} (\gamma), n\in 2\mathbb Z_++1$
to be a map from $\tRG$ to itself as the composition
\[
  R^-(\BG m)  \stackrel{ Res }{\longrightarrow}
   R^-(\BG n)\bigotimes R^- ( \BG{m -n})
 \stackrel{ \langle \sigma_n ( \g), \cdot
\rangle_{\wt}}{\longrightarrow}
 R^- ( \BG {m -n}).
\]
We denote by $\RGO $ the radical of the bilinear form $\langle \ ,
\ \rangle_{\wt}$ in $\tRG$ and denote by $\RGG$ the quotient $\tRG
/ \RGO$, which inherits the bilinear form $\langle \ , \
\rangle_{\wt}$ from $\tRG$. The following theorem is implied by
the identification of the Hopf algebra structures on $\tRG$ and
$\tSG$ (see Proposition~\ref{prop_hopf} below).

\begin{theorem}  \label{th_heis}
$\tRG$ is a module over the twisted Heisenberg algebra $\thg$ by
letting $ a_n (\g )$
 $( n \in 2\mathbb Z+1)$ act as $ \widetilde{a}_{ n} (\g )$
and $C$ as $1$. $\RGO$ is a submodule of $\tRG$ over $\thg$ and
the quotient $\RGG $
 is irreducible. The characteristic map
 $\ch$ is an isomorphism of $\tRG$ (resp. $\RGO$, $\RGG$)
 and $\tSG$ (resp. $\SGO$, $\tSGG$) as supermodules over $\thg$.
\end{theorem}
\subsection{The character $\chi_n (\g )$}

Recall that we have defined a map from $R(\G)$ to $R^-(\BG n )$ by
sending $\g$ to $\chi_n (\g) := \eta_n (\g) \circledast \chi_n$,
where $\chi_n$ is the character of the basic spin supermodule
$L_n$. The image of $\chi_n (\g)$ has the following elegant
description under the characteristic map in terms of a generating
function in a formal variable $z$.

\begin{proposition}  \label{prop_exp}
   For any $\g \in R(\G)$, we have
\begin{eqnarray}
 \sum\limits_{n \ge 0} \ch ( \chi_n( \g ) ) z^n
  &= & \exp \Biggl( \sum_{ n \ge 1, odd}
      \frac 2n \, a_{-n}(\g )z^n \Biggr). \label{eq_exp}
\end{eqnarray}
\end{proposition}
\begin{proof}
 Let $\g$ be a character of $\G$.It follows from Corollary
\ref{C:char} that

\begin{eqnarray*}
  \sum\limits_{n \ge 0} \ch ( \chi_n( \g ) ) z^n
  &= & \sum_{\rho \in {\mathcal OP}_n(\G_*)} 2^{l(\rho)}Z_{\rho}^{ -1}
         \prod_{c\in \G_*} \g (c)^{l (\rho(c))}
          a_{ -\rho (c) }' z^{|| \rho||}                  \\
  &= & \prod_{c\in \G_*} \Bigl ( \sum_{\lambda \in {\mathcal OP}_n}
         (2\zeta_c^{ -1}\g (c) )^{l (\lambda)}
         z_{\lambda}^{-1} a_{- \lambda} (c) z^{|\lambda|} \Bigr )
\\
  &= & \exp \Biggl  ( \sum\limits_{ n \geq 1, odd}
         \frac2n \sum\limits_{c \in \G_*}
           \zeta_c^{ -1} \g(c) a_{-n} (c) z^n \Biggl )      \\
  &= & \exp \Biggl( \sum_{ n \ge 1, odd}
         \frac 2n \, a_{-n}(\g )z^n \Biggr).
 \end{eqnarray*}

It follows from Eq.~(\ref{eq_virt}) that $\sum\limits_{n \ge 0}
\ch ( \chi_n( \g ) ) z^n$ is multiplicative on $\g$. Thus given
two characters $\beta, \g$ of $\G$, we have

\begin{eqnarray*}
\sum\limits_{n \ge 0} \ch ( \chi_n( \beta -\g ) ) z^n
 &=& \sum\limits_{n \ge 0} \ch ( \chi_n( \beta ) ) z^n
 \sum\limits_{n \ge 0} \ch ( \chi_n(- \g ) ) z^n \\
 &=& \exp \Biggl( \sum_{ n \ge 1, odd}
         \frac 2n \, a_{-n}(\beta )z^n \Biggr)
     \exp \Biggl( - \sum_{ n \ge 1, odd}
         \frac 2n \, a_{-n}(\g )z^n \Biggr)  \\
 &=& \exp \Biggl( \sum_{ n \ge 1, odd}
         \frac 2n \, a_{-n}(\beta -\g)z^n \Biggr).
\end{eqnarray*}
Therefore the proposition is proved.
 \end{proof}

\begin{corollary}  \label{cor_char}
   The formula (\ref{E:charvalue}) holds for any $\g \in R(\G)$.
 In particular $\chi_n (\wt)$ is self-dual if $\xi$ is self-dual.
\end{corollary}

Component-wise, we obtain
\begin{equation*}
  \ch (\chi_n (\g ) )=
\sum\limits_{\rho } \frac {2^{l(\rho)}}{z_\rho }\,
             a_{-\rho}(\g ),
\end{equation*}
where the sum is over all the partitions
$\rho$ of $n$ into odd integers.
\subsection{Isometry between $\tRG$ and $\tSG$}
    In Sect.~\ref{S:Hopf} we introduced
the  Hopf algebra structure
on $\tRG$. On the other hand it is well known that
there exists a natural Hopf
algebra structure on the symmetric algebra
$\tSG$ with the usual multiplication
and the comultiplication $\Delta$ given by
\begin{equation}\label{E:hopf3}
  \Delta ( a_n (\g ))
   = a_n (\g ) \otimes 1 + 1 \otimes a_n (\g ), \qquad n\in 2\mathbb
Z+1.
\end{equation}

\begin{proposition} \label{prop_hopf}
  The characteristic map $ \ch: \tRG \longrightarrow \tSG$
 is an isomorphism of Hopf algebras.
\end{proposition}

\begin{proof}
First the map $ch$ is a vector space isomorphism by comparing
dimension. The algebra isomorphism follows from the Frobenius
reciprocity. On the other hand, one can check directly that
 $$ \Delta (\sigma_n (\g)) =\sigma_n (\g ) \otimes 1
 + 1 \otimes \sigma_n (\g ), \qquad \g \in \G^*, n\in 2\mathbb Z+1. $$
Since $a_n (\g)$ (resp. $\sigma_n (\g)$) for $\g \in \G^*, n\in
2\mathbb Z+1$ generate $\tSG$ (resp. $\tRG$) as an algebra, we
conclude that $ch$ is a Hopf algebra isomorphism by
(\ref{E:hopf3}).
\end{proof}

Recall that we have defined a bilinear form $\langle \  ,  \,
\rangle_{\wt }$ on $\tRG$ and a bilinear form $\langle \ , \,
\rangle_{\wt }'$ on $\tSG$, and thus an induced one on $\tSG
\otimes \tSG$. The lemma below follows from our definition of
$\langle \ , \, \rangle_{\wt }'$ and the comultiplication
$\Delta$.

\begin{lemma}
The bilinear form $\langle \  ,  \, \rangle_{\wt } '$ on $\tSG$
can be characterized by the following two properties:

 1). $\langle a_{ -n} (\beta ), a_{ -m} (\g ) \rangle_{\wt}'
  = \frac{n}{2}\delta_{n, m} \langle \beta , \g  \rangle_{\wt}' ,
  \quad \beta, \g \in \G^*, m, n\in 2\mathbb Z+1.$

 2). $ \langle f g , h \rangle_{\wt}'
       = \langle f \otimes g, \Delta h \rangle_{\wt}' ,$
 where $f, g, h \in \tSG $.
\end{lemma}

\begin{theorem}  \label{th_isometry}
  The characteristic map $\ch$ is an isometry from the space
 $ (\tRG, \langle \ \ , \ \  \rangle_{\wt } )$ to
 $ (\tSG, \langle \ \ , \ \  \rangle_{\wt }' )$.
\end{theorem}

\begin{proof} Let $f$ and $g$ be any two super
class functions in $R^-(\BG n)$. By definition of the
characteristic map (\ref{E:ch}) it follows that
\begin{align*}
\langle ch(f), ch(g)\rangle_{\xi}' &=\sum_{\rho, \rho'\in
\mathcal{OP}_n (\G_*)} \frac 1{Z_{\rho}Z_{\rho'}}
f(\rho)g(\rho')\langle a_{-\rho}', a_{-\rho'}'\rangle_{\xi}'\\
&=\sum_{\rho, \rho' \in \mathcal{OP}_n (\G_*)} \frac
1{Z_{\rho}Z_{\rho'}}
f(\rho)g({\rho'})\frac{Z_{\rho}}{2^{l(\rho')}}
\prod_{c\in\G_*}\xi(c)^{l(\rho(c))} \delta_{\rho,
\overline{\rho'}}\\ &=\langle f, g\rangle_{\xi},
\end{align*}
where we have used the inner product identity
(\ref{eq_inner}) and (\ref{E:innerprod}).
\end{proof}

From now on the bilinear form $\langle \ , \ \rangle_{\wt}$ on
$\tRG$ is identified with the one $\langle \ , \ \rangle_{\wt} '$
on $\tSG$.
\section{Vertex representations and $\tRG$} \label{sec_vertex}
In this section, we construct twisted vertex operators on a space
$\bFG$ obtained by $\tRG$ tensored with a certain group algebra
constructed from the finite group $\G$. In a most important case
when $\G$ is a finite subgroup of $SL_2(\C)$, we obtain by vertex
operator techniques a twisted affine Lie algebra and toroidal Lie
algebra acting on $\bFG$ in terms of group-theoretic operators.
\subsection{A central extension of $\Rz/2\Rz$}
>From now on we assume that $\xi$ is a self-dual virtual character
of $\G$, then $\Rz$ is an integral lattice with the symmetric
bilinear form $\langle \ , \ \rangle_{\xi}$.

The quotient lattice $\Rtz=\Rz/2\Rz$ is an $(r+1)$-dimensional
vector space over the field $\mathbb F_2= \Z_2$. Write
$\overline{\alpha}=\alpha+2\Rz$. Then
 $$c_1(\oa, \ob)=\langle
\alpha ,\be \rangle_{\xi}+\langle\alpha ,\alpha\rangle_{\xi}
\langle\be ,\be\rangle_{\xi}\bmod{2}
  $$
is a natural (even) alternating form on $\Rtz$ and let $r_0$ be
its rank over $\mathbb F_2$. The alternating form $c_1$ defines a
central extension $\Rtzh$ of the abelian group $\Rtz$ by the
two-element group $\langle \pm 1\rangle$ (cf. \cite{FLM1}):
\begin{equation}
1\rightarrow \langle \pm 1\rangle\hookrightarrow \Rtzh
\stackrel{\Breve{}}{\rightarrow}\Rtz \rightarrow 1,
\end{equation}
such that $aba^{-1}b^{-1}=(-1)^{c_1(\breve{a}, \breve{b})}$, $a, b\in
\Rtzh$.
It is easily seen that
$\{\pm e_{\overline{\alpha}}\vert \alpha\in \Rz\}$ form a basis for
$\Rtzh$.
We note that $(e_{\oa})^2=1$ and $dim(\Rtzh)=2^{r+2}$.

Let $\Phi$ be a subgroup of $\Rz$ which is maximal
such that the alternating form $c_1$ vanishes on $\Phi/2\Rz$.

\begin{lemma} \cite{FLM2, FJW2} \label{L:centralext}
There are $2^{(r+1-r_0)}$ irreducible
$\Rtzh$-module structures on the
space $\mathbb C[\Rz/\Phi]$ such that $-1\in\Rtzh$ acts faithfully and

\begin{equation} \label{E:centralext}
e_{\oa}e_{\ob}=e_{\ob}e_{\oa}(-1)^{c_1(\oa, \ob)}
\end{equation}
as operators on $\mathbb C[\Rz/\Phi]$. Moreover, $\dim (\mathbb
C[\Rz/\Phi])=2^{\frac 12r_0}$.
\end{lemma}

We will denote the elements of $\mathbb C[\Rz/\Phi]$
by $e^{[\alpha]}$,
where $[\alpha]=\alpha+\Phi \in \Rz/\Phi$. Clearly

\begin{equation*}
e^{2[\alpha]}=1, \quad e^{[\alpha+\beta]}=e^{[\alpha]}e^{[\beta]}.
\end{equation*}
For $\alpha, \beta\in \Rz$ we write the action
of $\Rtzh$ on $\mathbb C[\Rz/\Phi]$ as
\begin{equation}\label{E:cocycle2}
e_{\oa}.e^{[\beta]}=\ep(\alpha, \beta)e^{[\alpha+\beta]}.
\end{equation}
Then one can check that $\ep$ is a well-defined
cocycle map
 from $\Rz\times \Rz\rightarrow \langle\pm 1\rangle$.
One also has $\ep(\alpha, \beta)=\ep(\alpha, -\beta)$.
\subsection{Twisted Vertex Operators $X ( \g, z)$}

We fix an irreducible $\Rtzh$-module structure on $\mathbb
C[\Rz/\Phi]$ described in Eqn. (\ref{E:cocycle2}). We extend the
actions of $e_{\oa}$ to the space of tensor product $$\bFG = \tRG
\bigotimes \mathbb C[\Rz/\Phi], $$ by letting them act on the
$\tRG$ part trivially.

Introduce the operators $ H_{ \pm n}( \g ), \g \in R(\G), n>0, $
as the following compositions of maps:
\begin{eqnarray*}
  H_{ -n} ( \g ) &:&
    R^- ( \BG m )
  \stackrel{\chi_n (\g) \otimes}{\longrightarrow}
    R^- ( \BG n\htimes\BG m )
  \stackrel{ {Ind} }{\longrightarrow}
    R^- ( \BG {n +m} )   \\
  H_n ( \g ) &:&
    R^- ( \BG m )
   \stackrel{ {Res} }{\longrightarrow}
    R^-( \BG n\htimes\BG {m -n})
   \stackrel{ \langle \chi_n (\g),
\cdot \rangle_{\wt} }{\longrightarrow}
    R^-( \BG{m -n}).
\end{eqnarray*}

We define
\begin{equation*}
  H_+ (\g, z) = \sum_{ n> 0} H_{ -n} ( \g ) z^n, \quad
  H_- (\g, z) = \sum_{ n> 0} H_{ n} (\g )z^{ -n} ,
\end{equation*}
where $z$ is a formal variable. We define the twisted vertex
operators $X_n (\g ), n \in {\mathbb Z}$, $\g\in\RG$ by the
following generating functions:
\begin{equation}  \label{eq_vo}
 X( \g, z)=H_+ (\g , z) H_- (\g , -z) e_{\overline{\g}}
=\sum\limits_{n \in
  {\mathbb Z}}
 X_n( \gamma) z^{ -n}.
\end{equation}
We note that
$X(-\ga, z)=X(\ga, -z)$.
The operators $X_n (\g )$
are well-defined operators acting on the space $ \bFG.$
We extend the bilinear form
$\langle \ , \ \rangle_{\wt}$ on $\tRG$ to $\bFG $ by letting
\[
  \langle f e^{[\alpha]}, g e^{[\beta]}\rangle_{ \wt} =
   \langle f, g \rangle_{\wt} \delta_{[\alpha],[\beta]}, \quad
    f, g\in \tRG, \alpha, \beta\in\Rz.
\]
We extend the $\mathbb Z_+$-gradation
>from $\tRG$ to
$\bFG$ by letting
\begin{eqnarray*}
  \deg  a_{ -n} (\g ) = n , \quad
  \deg e_{\overline{\g} } = 0.
\end{eqnarray*}

Similarly we extend the bilinear form
$\langle \  , \ \rangle_{\wt }$  to the space
$$
  \tVG = \tSG \bigotimes \mathbb C[\Rz/\Phi]
$$ and extend the $\mathbb Z_+$-gradation on $\tSG$ to a $\mathbb
Z_+$-gradation on $\tVG$.

The characteristic map $\ch$ will be extended to an isometry from
$\bFG$ to $\tVG$ by fixing the subspace $\mathbb C[\Rz/\Phi]$. We
will denote this map again by $\ch$.
\subsection{An identity for twisted vertex operators}
    We extend the characteristic map $ch$ to a linear map
$ch$: $End(\tRG)\rightarrow End(\tSG)$ by
\begin{equation}
ch(f).ch(v)=ch(f.v), f\in End(\tRG), v\in\tRG.
\end{equation}
The relation between the vertex operators defined
in (\ref{eq_vo}) and the Heisenberg algebra $\hg $ is revealed
in the following theorem.

\begin{theorem} \label{th_exp}
For any $\g \in R(\G)$, we have
  \begin{eqnarray*}
   \ch \bigl ( H_+ (\g, z) \bigl )
   &=& \exp \biggl ( \sum\limits_{ n \ge 1, \, odd} \frac 2n \,
    a_{-n} ( \g ) z^n \biggr ), \\
   \ch \bigl ( H_- (\g , z) \bigl )
   &=& \exp \biggl ( \sum\limits_{n \ge 1,\, odd}\frac 2n \,
    a_n (\g ) z^{-n}\biggr ).
  \end{eqnarray*}
\end{theorem}

\begin{proof} The first identity follows from
Proposition \ref{prop_exp} and the second one
is obtained by noting that
$H_+ (\g , z)$
 is the adjoint operator of $ H_- (\g , z^{ -1})$
with respect to the bilinear form
 $\langle \ , \ \rangle_{\wt}$.
\end{proof}

As a consequence we have
\begin{eqnarray*}
  && \ch \bigl ( X( \g , z)\bigl )  \\
  &= & \exp \biggl ( \sum\limits_{ n \ge 1, \, odd}
  \frac 2n \, a_{-n} ( \g ) z^n \biggr ) \,
  \exp \biggl ( -\sum\limits_{ n \ge 1, \, odd}
  \frac 2n \,{ a_n( \g)} z^{ -n} \biggr )e_{\overline{\g}}.
\end{eqnarray*}
Thus the characteristic map identifies the twisted vertex
operators $X(\g, z)$ defined via finite groups $\BG n$ with the
usual twisted vertex operators of \cite{FLM1, FLM2}.
\subsection{Product of two vertex operators}
The normal ordered product $: X(\alpha, z) X(\beta, w) :$,
$\alpha, \beta \in R(\G)$ of two vertex operators is defined
as follows:
\begin{equation*}
: X(\alpha, z) X(\beta, w) : =
H_+(\alpha, z)H_+(\beta, w) H_-(\alpha, -z)H_-(\beta, -w)
e_{\oa+\ob}.
\end{equation*}
The following theorem can be verified using the standard vertex
operator claculus (see e.g. \cite{FLM2, J, FJW2}), where the term
$\big(\frac{z -w}{z+w} \big)^{ \langle \alpha, \beta
\rangle_{\wt}}$ is understood as the power series expansion in the
variable $w/z$.
\begin{theorem}  \label{th_ope}
 For $\alpha, \beta \in R(\G)$ one has the following
operator product expansion identity for
twisted vertex operators.
 \begin{eqnarray*}
  X(\alpha, z) X(\beta, w) & =& \ep(\alpha, \beta)
  :X(\alpha, z) X(\beta, w):
  \biggr(\frac{z -w}{z+w}\biggr)^{ \langle \alpha, \beta
\rangle_{\wt}}.
 \end{eqnarray*}
\end{theorem}

The next proposition follows readily from Theorem~\ref{th_exp}.

\begin{proposition}  \label{prop_prim}
  Given $\alpha \in R(\G), \beta \in \Rz$ and $n\in 2\mathbb Z+1$, we
have
  $$
   [ a_n (\alpha ), X(\beta, z)]
   = \langle \alpha, \beta \rangle_{\wt} X(\beta, z) z^n.
  $$
\end{proposition}
\subsection{Twisted affine algebra $\tloopg$
and toroidal algebra $\thhg$}

Let $\mathfrak g$ be a rank $r$ complex simple Lie algebra of ADE
type, and let $\overline{\Delta}$ be the root system generated by
a set of simple roots $\alpha_1. \ldots, \alpha_r$. Let
$\alpha_{max}$ be the highest root. The Lie algebra is generated
by the Chevalley generators $e_{\alpha_i}, e_{-\alpha_i},
h_i=h_{\alpha_i} $. We normalize the invariant bilinear form on
$\mathfrak g$ by $(\alpha_{max}, \alpha_{max}) = 2$.

The twisted toroidal algebra $\thhg$ (associated to $\mathfrak g$)
is the associative algebra generated by (\cite{J2, FJW2}) $$
 C, h_i (m), x_n (\pm\alpha_i), m\in 2\mathbb Z+1,
n \in \mathbb Z, i =0, \ldots, r;
$$
subject to the relations: $C$ is central, $x_n(\alpha_i)= (-1)^n
x_n(-\alpha_i)$ and
\begin{align} \nonumber
   &{[ h_i (m), h_j (m')]}
  =\frac m2 a_{ij} \delta_{m, -m'}C,    \nonumber  \\
   &{[h_i (n), x_m (\alpha_j)]}
  =a_{ij} x_{n +m} (\alpha_j),  \nonumber    \\
   &{[x_n (\alpha_i), x_{n'}(- \alpha_i)]}
  = 8\{ h_i (n +n') + n\delta_{n, -n'} C \},
      \label{eq_pres}  \\
&{\sum_{s=0}^{a_{ij}}\binom{a_{ij}}{s}[x_{n+s}(\alpha_i),
x_{n'-a_{ij}-s}(\alpha_j)]}=0, \quad\mbox{if $a_{ij}\geq 0$}
\nonumber \\ &{\sum_{s=0}^{-a_{ij}}(-1)^s
\binom{-a_{ij}}{s}[x_{n+s}(\alpha_i), x_{n'-a_{ij}-s}(\alpha_j)]}=
0, \quad\mbox{if $a_{ij} <0$} \nonumber
\end{align}
where $n, n' \in \mathbb Z$, $m, m'\in 2\mathbb Z+1$, $i, j = 0,
1, \ldots, r$, and $h_i(2n)=0$ for $n\in\mathbb Z$.

The twisted affine algebra (cf. \cite{FLM1, FJW2}), denoted by
$\tloopg$, can be identified with the subalgebra of $\thhg$
generated by $$
 C, h_i (m), x_n (\pm\alpha_i), m\in 2\mathbb Z+1,
n \in \mathbb Z, i =1, \ldots, r; $$

The {\em basic twisted representation}  $V$ of $\tloopg$ is the
irreducible highest weight representation generated by a highest
weight vector which is annihilated by $ h_i (m) (m\in 2\mathbb Z_+
+1), x_n (\pm\alpha_i) (n \in \mathbb Z_+),$ and $C$ acts on $V$
as the identity operator.
%
%
%
%
%
\subsection{Realization of the twisted representations}
Let $\G$ be a finite subgroup of $SL_2(\mathbb C)$ and the virtual
character $\wt$ to be twice the trivial character minus the
character of the two-dimensional defining representation of
$\G\hookrightarrow SL_2(\mathbb C)$. The following is the
well-known list of finite groups of $SL_2(\mathbb C)$: the cyclic,
binary dihedral, tetrahedral, octahedral and icosahedral groups.
McKay observed that the associated matrix to the weighted bilinear
form $\langle-,-\rangle_{\wt}$ on $R(\G)$ with respect to the
basis of irreducible characters can be identified with an affine
Dynkin diagram of ADE type \cite{Mc}.

The following theorem provides a finite group realization of the
vertex representation of the twisted toroidal Lie algebra $\thhg$
on $\bFG$.

\begin{theorem}  \label{th_mck}
 A vertex representation of the twisted toroidal Lie algebra $\thhg$
 is defined on the space $\bFG$ by letting
 \begin{align*}
  x_n (\alpha_i) \mapsto X_n (\g_i), &\qquad
  x_n (-\alpha_i) \mapsto  \ep(\g_i, \g_i)X_n (-\g_i),\\
  h_i (m) \mapsto a_m (\g_i), &\qquad C\mapsto 1,
 \end{align*}
 where $n \in \mathbb Z, m\in 2\mathbb Z+1, 0 \leq i \leq r$.
\end{theorem}

\begin{proof} All the commutation relations without binomial
coefficients are easy consequences
of Proposition~\ref{prop_prim} and Theorem~\ref{th_ope}
by the usual vertex operator calculus in the twisted picture
 (see \cite{FLM2, J}). The corresponding relations with binomial
coefficients
in $\tVG$ are equivalent to
\begin{align*}
(z+w)^{a_{ij}}[X(\g_i, z), X(\g_j, w)]&=0, \qquad
\mbox{$a_{ij}\geq 0$},\\ (z-w)^{-a_{ij}}[X(\g_i, z), X(\g_j,
w)]&=0, \qquad \mbox{$a_{ij} > 0$}.
\end{align*}
This is proved by using Theorem~\ref{th_ope} by a standard method
in the theory of vertex algebras (cf. e.g. \cite{FLM2, J2}).
\end{proof}

Note that the quotient lattice $\Rz / R^0_{\mathbb Z}$ inherits a
positive definite integral bilinear form from that of $\Rz$. Let
$\Gbar= \{ \g_1, \g_2, \ldots, \g_r \}$ be the set of non-trivial
irreducible characters of $\G$. Let $R_{\mathbb Z}(\Gbar)$ be the
sublattice of $\Rz$ generated by $\Gbar$. Denote by $Sym(\Gbar)$
the symmetric algebra generated by $a_{-n} (\g_i),$ $ n\in
2\mathbb Z_++1, $ $i =1, \ldots , r$, which is isometric to
$\tSGG$ and $\RGG$ as well. The irreducible $\Rtzh$-module
$\mathbb C[\Rz/\Phi]$ induces an irreducible $\Rtzhbar$-module
structure on $\mathbb C[R_{\mathbb Z}(\Gbar)/\overline{\Phi}]$
given by the restriction of $c_1$. Denote by $\overline{r}_0$ its
rank. In this case if the determinant of the Cartan matrix is an
odd integer (see Lemma \ref{L:centralext}), then
$\overline{r}_0=0$ and the space $\mathbb C[R_{\mathbb
Z}(\Gbar)/\overline{\Phi}]$ is trivial.

We define
\begin{align*}
 \tVGG  &= \tSGG \bigotimes \mathbb C[R_{\mathbb
Z}(\Gbar)/\overline{\Phi}]
 \cong Sym (\Gbar) \bigotimes \mathbb C[R_{\mathbb
Z}(\Gbar)/\overline{\Phi}] ,\\ \bFGG  &= \tRGG \bigotimes \mathbb
C[R_{\mathbb Z}(\Gbar)/\overline{\Phi}].
\end{align*}
Obviously $ch$ restricted to $\bFGG$ is an isometric
isomorphism onto
$\tVGG$.
We remark that $\bFG$ is isomorphic to the tensor
product of the space $\FGG$
associated to $\Rzz$ and the space associated to
the rank $1$ lattice $\mathbb Z \delta$
equipped with the zero bilinear form.

The identity for a product of vertex operators $X (\g, z)$
associated to $\g \in \overline{\Delta}$
(cf.~Theorem~\ref{th_mck}) implies that $\tVGG$ provides a
realization of the vertex representation of $\tloopg$ on $\tVGG$
(cf. \cite{FLM1}). The following theorem establishes a direct link
>from the finite group $\G \in SL_2(\mathbb C)$ to the affine Lie
algebra $\tloopg$.

\begin{theorem}
The operators $X_n(\g ), \g \in \overline{\Delta}, a_n (\g_i), i
=1, 2, \ldots, r, n \in \mathbb Z$ define an irreducible
representation of the affine Lie algebra $\tloopg$ on $\bFGG$
which is isomorphic to the twisted basic representation.
\end{theorem}
\section{The character tables of $\BG n $ and vertex operators}
\label{sec_char}
     In this section we outline how to use
the specialization of $\wt=\g_0$ to obtain the character table for
the spin supermodules of $\BG n$ from our vertex operator
approach. The proofs are similar to those given in
\cite{FJW2}(also cf. \cite{J}) which we omit.
\subsection{Components of vertex operators}
Set $\wt=\g_0$. The weighted bilinear form reduces to the standard
one and $\Rz\simeq\mathbb Z^{r+1}$.

Let $A=(1-\delta_{ij})_{i,j=0}^r$, the matrix of the alternating
form $c_1$ over $\mathbb F_2$, then $A^2=\overline rI$. Here
$\overline r=0 $ if $r$ is even and $1$ if $r$ is odd.
Consequently Lemma \ref{L:centralext} implies there are exactly
$2^{\overline{r+1}}$ irreducible $\hat{R}^-_{\mathbb
F_2}(\G)$-module structures on the
$2^{\lceil\frac{r+1}2\rceil}$-dimensional space $\mathbb
C[R_{\mathbb Z}(\G)/\Phi]$, where $\lceil\frac{r+1}2\rceil$
denotes the smallest integer $\leq\frac{r+1}2$.
One of the (at most) two irreducible
module structures is given by the cocycle $\ep(\g_i, \g_j)=1$, for
$i\leq j$, and $\ep(\g_i, \g_j)=-1$, for $i>j$.

In the following result
the bracket
$\{ \ , \ \}$ denotes the anti-commutator.

\begin{theorem}  \label{th_cliff}
  The operators $X^+_n (\g_i), X^-_n (\g_i)$
 $(n \in \mathbb Z, 0 \leq i \leq r)$
 generate a Clifford algebra:
 \begin{eqnarray}
\{X_n( \g_i), X_{n'}(\g_j)\} &=&2(-1)^n\delta_{ij}\delta_{n, -n'},
\nonumber   \\ \{X_n(-\g_i), X_{n'}(-\g_j)\}
&=&2(-1)^n\delta_{ij}\delta_{n, -n'}, \label{eq_delta}\\
\{X_n(\g_i), X_{n'}(-\g_j)\} &= & 2\delta_{ij}\delta_{n, -n'}.
  \nonumber
 \end{eqnarray}
\end{theorem}
\begin{proof}
The relations can be established using Theorem \ref{th_ope} by
means of standard techniques in the theory of vertex algebras, see
\cite{J} for detail for a similar circumstance.
\end{proof}
\subsection{Spin character tables of $\BG n$ and vertex
operators}
\label{S:chartab}
Let $R^-_{\mathbb Z}(\BG n)$ be the lattice
generated by the characters of spin irreducible
$\BG n$-supermodules. Then
$R^-_{\mathbb Z}(\BG n)\otimes\mathbb C\simeq R^-(\BG n)$.

For an $m$-tuple index $\phi=(\phi_1, \cdots, \phi_m)\in
\mathbb Z^m$ we denote
\begin{equation*}
X_{\phi}(\g )=X_{\phi_1}(\g )\cdots X_{\phi_m}(\g )e^{[\alpha]},\quad
x_{\phi}(\g)=(X_{\phi_1}(\g ).1)\cdots (X_{\phi_m}(\g ).1).
\end{equation*}

The space $\tRG$ is spanned by $X_{\phi}(\g )$. However
$X_{-n}(\pm\g ).e^{[\alpha]}=0$ if $n<0$. Using the Clifford algebra
structure we know that $\tRG$ is spanned by
$X_{\phi, [\alpha]}(\g )$, where
$\phi$ runs through $m$-tuples $\phi=(\phi_1, \cdots, \phi_m)\in
\mathbb Z^m$ and $[\alpha]\in R_{\mathbb Z}(\G)/\Phi$. We define the raising operator $R_{ij}$ by
\begin{equation*}
R_{ij}(\phi_1, \cdots, \phi_m)=(\phi_1, \cdots, \phi_i+1,
\cdots, \phi_j-1, \cdots, \phi_m).
\end{equation*}
Then we define the action of the raising operator $R_{ij}$
on $X_{\phi, [\alpha]}(\g )$ or $x_{\phi}(\g)$ by
$X_{R_{ij}\phi, [\alpha]}(\g )$ or $x_{R_{ij}\phi}(\g)$ respectively.

Given $ \la \in \mathcal {OP}(\G^*)$, we define
\[
X_{\lambda}
  = \prod_{i=0}^r X_{-\la( \g_i )} (\g_i).
\]
Similarly we write
$x_{\lambda}=\prod_{\g \in \mathcal \G^*}x_{\la( \g )} (\g).$

\begin{proposition} \cite{FJW2} \label{T:char1}
 The vectors $X_{\lambda}e^{[\alpha]}$ for
 $ \lambda=(\la(\g ))_{\g \in \G^*} \in \mathcal {SP} (\G^*)$
and $[\alpha]\in R_{\mathbb Z}(\G)/\Phi$
 form an orthogonal basis in the vector space $\tRG$ with $\langle
X_{\lambda}e^{[\alpha]}, X_{\mu}e^{[\beta]}\rangle=
2^{l(\lambda)}\delta_{\lambda, \mu}\delta_{[\alpha], [\beta]}$.
Moreover, we have that

\begin{equation} \label{E:raisingop}
X_{\lambda}e^{[\alpha]}=\prod_{\g\in\G_*}
\prod_{ij}\frac{1-R_{ij}}{1+R_{ij}}x_{\la}e^{l(\lambda)[\alpha]}
=(x_{\lambda}+\sum_{\lambda\gg\mu}c_{\la, \mu}x_{\mu})e^{l(\lambda)[\alpha]},
\end{equation}
where $c_{\la, \mu}\in \mathbb Z$, and $R_{ij}$ is the raising
operator.
\end{proposition}

We remark that the basis elements in $\tSG$ corresponding to
$x_{\la}$ are the classical symmetric functions $Q_{\la}$ called
Schur's Q-functions \cite{M, J}. If we take $a_{ -n } (\g ), n >0
, \g \in \G^*$ as the $n$-th power sum in a sequence of variables
$ y_{ \g } = ( y_{i\g } )_{i \geq 1}$, then the space $\tSG$
becomes a distinguished subspace of symmetric functions generated
by odd degree power sums indexed by $ \G^*$. In particular
$Q_{\lambda} (\g)$ is the Schur's Q-function in the variables
$y_{\g}$. For $\lambda \in \mathcal P (\G^*)$, we denote
\begin{eqnarray}  \label{eq_schur}
  Q_{\lambda} = \prod_{\g \in \G^*} Q_{\lambda(\g)} (\g) \in \tSG.
\end{eqnarray}

For $\la\in\mathcal {SP}_n( \G^*)$, we define
$\BG {\la}=\BG{\la(\ga_0)}\htimes \cdots\htimes \BG{\la(\ga_r)}$.
For a partition $\mu$ and an irreducible character $\ga$ of
$\G$, we define the
spin character $\chi_{\mu}(\g)$ of $\BG{\mu}$ to be
$\chi_{\mu_1}(\ga)\otimes \cdots\otimes
\chi_{\mu_l}(\ga)$ (see Corollary \ref{C:char}).

\begin{theorem}\label{T:main}
For each strict partition-valued function $\lambda \in \mathcal
{SP}_n(\G^*)$, the vector $2^{-[ l(\la)/2]}Q_{\lambda}$
corresponds, under the characteristic map $ch$,  to the
irreducible character $\chi_{\la}$ of the spin $\BG n$-supermodule
given by a $\mathbb Z$-linear combination of

\begin{equation}\label{E:main3}
Ind_{\BG{\rho}}^{\BG n}
\chi_{\rho(\g_0)}(\g_0)\otimes \cdots\otimes\chi_{\rho(\g_{r})}(\g_r),
\end{equation}
where $\rho\ll\la$ and the first summand is $\rho=\la$ with
multiplicity one. Its character at the conjugacy class of type
$(\mu, \emptyset)$, where $\mu \in \mathcal {OP}_n(\G_*)$, is
equal to
 the matrix coefficient
 \begin{eqnarray}\label{E:table}
  2^{l(\mu)-\lceil l(\la)/2\rceil}\langle
    X_{\lambda}e^{-[\la]}, a_{-\mu}'
  \rangle,
 \end{eqnarray}
where $[\la]=\sum_{i=0}^rl(\la(\ga_i))[\ga_i]$. Moreover the
degree of the character is equal to

\begin{equation}\label{E:degree}
2^{n-\lceil l(\la))/2\rceil}n!\prod_{\g\in \G^*}
\big(\frac{deg(\g)^{|\la(\g)|}} {\prod_{1\leq i\leq l(\la(\g))}
\la_i(\g)!}
\prod_{ij}\frac{\la_i(\g)-\la_j(\g)}{\la_i(\g)+\la_j(\g)} \big).
\end{equation}
\end{theorem}

All the irreducible spin $\BG n$-supermodules can be described
easily as follows.  For each irreducible character $\g\in\G^*$ let
$U_{\g}$ be the irreducible $\G$-module affording $\g$. For each
strict partition $\nu$ of $n$ let $V_{\nu}$ be the corresponding
irreducible spin supermodule of $\B n$ (cf. \cite{Sg, Jo2}). We
have seen earlier that $U^{\otimes n}_{\g}\otimes V_{\nu}$ is an
irreducible spin $\BG n$-supermodule.

\begin{proposition}
For each strict partition-valued function $\la=(\la(\g))$ $\in$
$\mathcal{SP}_n(\G^*)$, the tensor product $$ \prod_{\g\in\G^*}
\big(U^{\otimes l(\la(\g))}_{\g}\otimes V_{\la(\g)}\big) $$
decomposes completely into $2^{\lceil m/2\rceil}$ copies of an
irreducible spin $\BG{\la}$-super\-module, where $m$ denotes the
number of the partitions of odd length among $\la(\g)$. Denote
this irreducible module by $W_{\la}$. Then the induced supermodule
$Ind_{\BG{\la}}^{\BG n}W_{\la}$ is the irreducible spin $\BG
n$-supermodule corresponding to $\la$, and it is of type $M$ or
$Q$ according to $l(\la)$ is even or odd.
\end{proposition}

\bibliographystyle{amsalpha}

\end{document}